\documentclass{cmslatex}
\usepackage[paperwidth=7in, paperheight=10in, margin=.875in]{geometry}
 \usepackage[backref,colorlinks,linkcolor=red,anchorcolor=green,citecolor=blue]{hyperref}
\usepackage{amsfonts,amssymb}
\usepackage{amsmath}
\usepackage{graphicx}
\usepackage{cite}
\usepackage{enumerate}

\usepackage{color,xcolor}
\usepackage{graphics,subfigure}
\usepackage{url}
\usepackage{hyperref}

\usepackage{algorithmic,algorithm}
\renewcommand{\theequation}{\arabic{section}.\arabic{equation}}

\newcommand{\dpar}[2]{\dfrac{\partial #1}{\partial #2}}

\newcommand{\R}{\mathbb R}
\newcommand{\Z}{\mathbb Z}
\newcommand{\N}{\mathbb N}

\renewcommand{\P}{\mathbb P}

\newcommand{\bbf}{{\mathbf {f}}}

\newcommand{\bu}{\mathbf{u}}
\newcommand{\bv}{\mathbf{v}}

\newcommand{\BF}{\mathbf{F}}
\newcommand{\BBF}{\mathbf{\mathcal{F}}}

\newcommand{\bbu}{\mathbf{u}}

\newcommand{\hf}{\hat{\mathbf{f}}}

\newcommand{\ba}{\mathbf{a}}

\newcommand{\hu}{\hat{u}}
\newcommand{\hv}{\hat{v}}
\newcommand{\M}{\mathbb M}
\newcommand{\revA}[1]{{\color{black}#1}}
\newcommand{\revB}[1]{{\color{black}#1}}
\newcommand{\revD}[1]{{\color{black}#1}}
   \allowdisplaybreaks
\begin{document}
 \title{Some preliminary results on a high order asymptotic preserving \revA{ computationally explicit} kinetic scheme.\thanks{Received date, and accepted date (The correct dates will be entered by the editor).}}


          \author{R\'emi Abgrall\thanks{Institute of Mathematics and Institute of Computational Sciences,	Universit\"at Z\"urich, Wintherturerstrasse 190, Z\"urich, Switzerland, (\href{mailto:remi.abgrall@math.uzh.ch}{remi.abgrall@math.uzh.ch}).} 
          \and Davide Torlo \thanks{Inria Bordeaux - Sud-Ouest, 200 avenue de la vieille tour, 33405 Talence, France, (\href{mailto:davide.torlo@inria.fr}{davide.torlo@inria.fr}), \url{davidetorlo.it}.}}

         \pagestyle{myheadings} \markboth{Results on Asymptitic Preserving Kinetic Scheme}{R. Abgrall and D. Torlo} \maketitle

          \begin{abstract}
 	In this short paper, we intend to describe one way to construct arbitrarily high order  kinetic schemes on regular meshes. The method can be arbitrarily high order in space and time,  run at least CFL one, is asymptotic preserving \revA{and computationally explicit, i.e., the computational costs are of the same order of a fully explicit scheme}. We also introduce a non linear stability method that enables to simulate  problems with discontinuities, and it does not kill the accuracy for smooth regular solutions.
          \end{abstract}

\section{Introduction}

Let us specify first the context. We are given the PDE
\begin{subequations}\label{eq:1}
	\begin{equation}
		\label{eq:1:1}
		\dpar{\bbu}{t}+\dpar{\bbf(\bbu)}{x}=0
	\end{equation}
	with the initial condition \begin{equation}\label{eq:1:2}
		\bbu(x, 0)=\bbu_0(x),
	\end{equation}
\end{subequations}
with $\bbu\in \R^p$ and $\bbf:\R^p\rightarrow \R^p$. a Lipschitz continuous flux
It is known, at least since the work of Jin \cite{Jin} and then Natalini \cite{Natalini} and co-workers, that this system can formally be seen as the limit for $\varepsilon\rightarrow 0$ of a relaxation system:
\begin{subequations}\label{eq:2}
	\begin{equation}\label{eq:2:1}
		\dpar{\BF}{t}+ \Lambda \dpar{\BF}{x}=\dfrac{\M(\P\BF)-\BF}{\varepsilon}
	\end{equation}
\end{subequations}
with $\BF\in \R^{k\times p}$, $\M$ is a Maxwellian and 
$\P$ is a linear operator such that $\P\M(\P\BF)=\P\BF$. The constant matrix $\Lambda$ and the flux $\bbf$ are linked by  $\P\Lambda\M(\P\BF)=\bbf(\P\BF)$. The simplest example, due to Jin and Xin \cite{Jin},  is 
\begin{equation*}
	\begin{array}{l}
		\dpar{u}{t}+\dpar{v}{x}=0\\
		\dpar{v}{t}+a^2 \dpar{u}{x}=\dfrac{f(u)-v}{\varepsilon}$$
	\end{array}
\end{equation*}
that can be rewritten in the form \eqref{eq:2} with:
\begin{equation}
	\begin{split}
		\dpar{f_1}{t}+a\dpar{f_1}{x}&=\dfrac{\M_1-f_1}{\varepsilon},\\
		\dpar{f_2}{t}-a\dpar{f_2}{x}&=\dfrac{\M_2-f_2}{\varepsilon},\\
	\end{split}
\end{equation}
i.e.,
where 
$$\BF=\begin{pmatrix}f_1\\ f_2\end{pmatrix},\quad  \Lambda=\begin{pmatrix}a & 0\\0&-a\end{pmatrix}, \quad \P\BF=f_1+f_2 \text{ and }\M=\begin{pmatrix} \M_1 \\ \M_2\end{pmatrix}$$ where the Maxwellian is defined from the relations
$$\M_1+\M_2=f_1+f_2=u, \qquad  a(\M_1-\M_2)=f(u), $$
i.e.
$$\M_1(f,a)=\frac{1}{2}\left (f_1+f_2+\dfrac{f(u)}{a}\right ), \M_2(f,a)=\frac{1}{2}\left (f_1+f_2-\dfrac{f(u)}{a}\right ).$$
We know that $a$ must be larger than the max of $|f'(u)|$ because of the Whitham sub-characteristic condition, obtained  via  a formal Chapman Enskog expansion. Another argument is, as shown by \cite{Bouchut}, that under this condition the two Maxwellian $\M_1$ and $\M_2$ satisfy a monotonicity condition, i.e. the BGK model becomes compatible with entropy inequalities.

\bigskip
The questions we address in this paper are the following: given a system \eqref{eq:1} and a regular grid of spatial step $\Delta x>0$, can we construct a \revA{computationally} explicit scheme that solves \eqref{eq:2}  with uniform accuracy of order $r>0$ for all $\varepsilon>0$ and with a CFL condition,  based on the matrix $\Lambda$, that is larger than $1$\footnote{ Initially, the first author was motivated by understanding in a better way the LBM method, even though the answer is not about the LBM method at all. The only remaining property between what we look for and the LBM method is the CFL condition.}.
The answer is yes, and this paper proposes a simple construction in one dimension. 
\revA{With \textit{computationally explicit} we mean that the solution of a certain scheme does not require any nonlinear solver, nor the inversion of a mass matrix.}

High order accurate methods for kinetic problems \`a la Shi-Jin  has received a lot of attention in the recent years. For a long time the state of the art was that of second  
second order in time and space  finite volume with TVD like stabilisation, see e.g \cite{AregbaNatalini}. For higher than second order,  one may mention \cite{MR3635826} where a splitting approach is adopted with a regular CFL stability condition for the overall finite volume scheme, \cite{MR1910754} where  relaxed upwind schemes are proposed running up to CFL $=1$ and up to third order in time/space, again in a finite volume context. In \cite{sead}, a WENO approach is proposed. In  \cite{COULETTE2019485} a discontinuous Galerkin approximation of the system \eqref{eq:2} is developed (with a temporal scheme allowing very large CFL number). 
In the kinetic literature where the fluid system is represented with the BGK approximation, so that dense and less dense flows can be simulated, there has also been a large effort towards high order schemes with asymptotic preserving properties. One may mention \cite{MR3033046} for hyperbolic systems with diffusion,  \cite{MR3969008} where a high order conservative semi-Lagrangian technique is developed. 

We want to go beyond that, with very simple and cheap numerical schemes that are potentially arbitrary high order and run at CFL $=1$, with an accuracy that is independent of the relaxation parameter $\varepsilon$.
The format of the paper is as follows. We first introduce the general method which amounts to describe the discretisation of $\Lambda \dpar{\BF}{x}$ and a time discretisation. We take into account the source term. The scheme resulting from this discretisation is fully implicit. The next step is to show that, thanks to the operator $\P$, and using  a particular time discretisation, we can make it  \revA{computationally} explicit, and high order accurate, independently of the parameter $\varepsilon$. Several choices of $\Lambda$ and Maxwellians $\M$ are described. We also address the question of the non linear stabilisation of the method when discontinuities appear. Several numerical examples, covering scalar and system cases, are then proposed to show the relevance of the method. The accuracy is checked for the scalar case.
\section{General discretisation principle}
Starting from \eqref{eq:2}, the idea is to discretise first in space $\Lambda \dpar{\BF}{x}$. This introduces an error which we assume to be $O(\Delta x^q)$,
\begin{equation}
	\label{model}
	\dpar{\BF}{t}+\dfrac{1}{\Delta x}\Lambda\delta\BF=\dfrac{\M(\P\BF)-\BF}{\varepsilon}+O(\Delta x^q).
\end{equation}
The second step is to discretise in time, so that we expect that the resulting scheme will be of order $p$ in space and time, at least for moderate values of $\varepsilon$. 
The problem is then two-fold: (i) how to define the discretisation operator $\delta$ for which a minimum requirement is the semi discrete linear stability when there is no source term, (ii) how to discretise in time so that the accuracy is uniform in time and  $\varepsilon$. We first discuss the issue of time discretisation, then space discretisation.

\subsection{Time discretisation}
One may use IMEX Runge-Kutta schemes, and more precisely SSP IMEX Runge-Kutta schemes, to have a better control of the stability properties of the method. 
Rewriting \eqref{eq:2:1} as the sum of a non stiff term and a stiff one
\begin{equation}
	\label{zebi}\dfrac{dU}{dt}+\mathcal{F}(U)=\frac{\mathcal{G}(U)}{\varepsilon}\end{equation}
an IMEX method is defined by two Butcher's tableaux
$$
\begin{array}{c|c}
	c & A\\ \hline
	0&b^T
\end{array} \quad \text {and } \quad \begin{array}{c|c}
	\tilde{c} & \tilde{A}\\ \hline
	0&\tilde{b}^T
\end{array}
$$
where the first one is for non stiff part, while the second one is for the stiff part:
\begin{equation}\label{Imex}
	\begin{split}
		U_0&=U^n\\
		\vdots& \\
		U_k&=U_0+\Delta t \sum_{j=1}^{k-1} a_{kj}\mathcal{F}(U_j)+\dfrac{\Delta t}{\varepsilon}\sum_{j=1}^{s} \tilde{a}_{kj}\mathcal{G}(U_j)\\
		\vdots\\
		U^{n+1}&=U^n+\Delta t\sum_{j=1}^s b_j\mathcal{F}(U_j)+\dfrac{\Delta t}{\varepsilon}\sum_{j=1}^{s} \tilde{b}_{j}\mathcal{G}(U_j)
	\end{split}
\end{equation}
with various compatibility conditions so that a given order is reached, see \cite[Chapter IV]{Hairer}. Anticipating a bit, if there exists a linear operator $\P$ such that $\P \mathcal{G}=0$ as here, we see that, applying $\P$ to \eqref{Imex}, a necessary condition is that the explicit RK scheme defined by the explicit part is itself SSP. Since we want to have a running CFL number of at least one, this needs that the SSP RK scheme must have a CFL number of at least $1+\epsilon$, $\epsilon>0$. \revD{To our knowledge there are some explicit SSP RK schemes satisfying this condition, inter alia \cite{Ruuth}, but they are not generalizable to arbitrarily high order of accuracy, and no IMEX versions are available.} 

For this reason, we use \revA{an IMEX deferred correction (DeC) method. \revB{It is a general way of building arbitrarily high order Runge Kutta schemes. It also allows more freedom in the spatial discretization, for instance, allowing the use of lumped mass matrix \cite{Abgrall}. Its implicit and IMEX versions allow to use a combination of more traditional low order IMEX schemes and arbitrarily high order implicit RK schemes, obtaining arbitrarily high order IMEX schemes}. We leave the study of SSP version of these schemes for future research. The final IMEX DeC scheme we obtain is computationally explicit and it is also matrix-free.}

\subsubsection{Deferred Correction}
\revB{The DeC is an iterative procedure that was proposed and developed in its explicit version in \cite{Dutt_DeC} and in an implicit version in \cite{minion}. It was applied to hyperbolic PDE, for instance, in \cite{Abgrall}, with a new formalism which makes the proof of its properties more straightforward. An IMEX version of this algorithm applied to hyperbolic PDE is available in \cite{Torlo}, and the algorithm we discuss in the following is a modification of this one. 
	With the notation of \cite{Abgrall}, the DeC uses two operators: one high order accurate $\mathcal{L}^2$, which defines a fully implicit method, and a low order easy to solve $\mathcal{L}^1$ operator.  The process allows to approximate with arbitrary accuracy the solution of the high order operator $\mathcal{L}^2$, with the simplicity of the operator $\mathcal{L}^1$.
	We start with the description of the high order operator $\mathcal{L}^2$.}

Let us consider $q+1$ points in $[0,1]$, $c_0=0<c_1\ldots<c_i<\ldots <c_q=1$ and the quadrature formula
$$\int_{t_n}^{t_n+c_i\Delta t}\varphi(s)\; ds\approx \Delta t \sum_{j=0}^q a_{ij}\varphi (t_n+c_j\Delta t).$$
More precisely, if $\{\ell_j\}$ are the Lagrange polynomials associated to the partition $\{c_j\}_{j=0}^q$, if we take
$$a_{ij}=\int_0^{c_i}\ell_j(s)ds,$$
the quadrature formula is of order $q+1$. 
We will always require that the quadrature formula are consistent, i.e.
\begin{equation}
	\label{consistant}
	\sum_{j=0}^qa_{ij}=c_i.
\end{equation} Considering $x_k$, a grid point, and setting $\BF_k^{n,j}\approx\BF(x_k,t_n+c_j\Delta t)$ and $\BF_k^{n,0}=\BF_k^n$, an approximation of \eqref{eq:2} is:
\begin{equation}
	\label{eq:3}
	\BF_k^{n,j}-\BF_k^{n,0}+\dfrac{\Delta t }{\Delta x}\bigg (\sum_{l=0}^q a_{il}\Lambda \delta_k \BF^{n,l}\bigg )-\mu \sum_{l=0}^q a_{il}\big ( \M\P\BF_k^{n,l}-\BF_k^{n,l}\big )=0, \quad j=1, \ldots , q
\end{equation}
where $\mu=\tfrac{\Delta t}{\varepsilon}$ and $\tfrac{\delta \BF}{\Delta x}$ is a consistent approximation of $\dpar{\BF}{x}$.  We will set $\BF_k^{n+1}=\BF_k^{n,q}$.
The relations \eqref{eq:3} can be rewritten in matrix form, setting
$$\BBF_k=\big (\BF_k^{n,1}, \ldots , \BF_k^{n,q}\big )^T, \quad \BBF_k^{(0)}=\big (\BF_k^{n,0}, \ldots , \BF_k^{n,0}\big )^T=\big (\BF_k^n, \ldots \BF_k^n\big ),$$ and neglecting the index of the timestep $n$, as
\begin{equation}\label{eq:3bis}
	\begin{split}
	\BBF_k\!\!-\BBF_k^{(0)}+\dfrac{\Delta t}{\Delta x} \Lambda A \delta_k\BBF\!&-\mu A\big ( \M(\P\BBF_k)\!-\BBF_k\big )\!+\!\dfrac{\Delta t}{\Delta x} \Lambda \ba_0 \otimes \delta_k\BF^{n,0}\!\!\\
		&-\mu\ba_0 \otimes (\M(\P\BF_k^{n,0})\!-\BF_k^{n,0})\!=0,
	\end{split}
\end{equation}
where by abuse of language we  have written 
$$\M(\P\BBF)=\big ( \M(\P\BF^1), \ldots , \M(\P\BF^q)\big )^T.$$ The matrix $A$ 
is
$$A=\begin{pmatrix}
	a_{11} &\ldots & a_{1q}\\
	\vdots & \vdots &\vdots\\
	a_{q1} & \ldots & a_{qq}
\end{pmatrix}$$
and we have $$\mathbf{a}_0=\begin{pmatrix} a_{0q}\\ \vdots \\ a_{01} \end{pmatrix}.$$

As a result, \eqref{eq:3bis} is implicit, and in general non linear, because of the Maxwellian.  In order to simplify the resolution, we consider a simpler scheme, where the source term discretisation  remains the same and the forward Euler method is used on each sub-time step:
\begin{equation}
	\label{eq:4}
	\BF_k^{n,j}-\BF_k^n+c_j\dfrac{\Delta t }{\Delta x}\Lambda \delta_k \BF^{n,0}-\mu \sum_{l=0}^q a_{jl}\big ( \M(\P\BF_k^{n,l})-\BF_k^{n,l}\big )=0, \quad j=1, \ldots , q.
\end{equation}
We rewrite this as:
\begin{equation}\label{eq:4bis}
	\BBF_k-\BBF_k^{(0)}+\dfrac{\Delta t }{\Delta x} C\Lambda\delta_k\BBF^{(0)}-\mu A\big ( \M(\P\BBF_k)-\BBF_k\big )-\mu\ba_0
	\otimes (\M\P\BF_k^{n,0}-\BF_k^{n,0})
	=0
\end{equation}
where $C=\text{ diag }\big ( c_1, \ldots, c_q\big )$ and $\BBF^{(0)}=(\BF^{n,0}, \ldots , \BF^{n,0})^T$.

This leads to the introduction of two operators $\mathcal{L}^1$ and $\mathcal{L}^2$ acting on $\BBF=(\ldots, \BBF_k, \BBF_{k+1}, \ldots)$ and defined as:
$$\big [\mathcal{L}^1(\BBF)\big ]_k:=\BBF_k-\BBF_k^{(0)}+\dfrac{\Delta t }{\Delta x} C\Lambda\delta_k\BBF^{(0)}-\mu A\big ( \M(\P\BBF_k)-\BBF_k\big )-\mu\ba_0 \otimes (\M(\P\BF^{n,0}_k)-\BF^{n,0}_k),$$
and
\begin{equation}
	\begin{split}
\big [\mathcal{L}^2(\BBF)\big ]_k:=&\BBF_k-\BBF_k^{(0)}+\dfrac{\Delta t}{\Delta x} \Lambda A \delta_k\BBF-\mu A\big ( \M(\P\BBF_k)-\BBF_k\big )+\dfrac{\Delta t}{\Delta x} \Lambda \ba_0 \otimes \delta_k\BF^{n,0} \\ &- \mu\ba_0  \otimes (\M(\P\BF^{n,0}_k)-\BF^{n,0}_k).
	\end{split}
\end{equation}
So that \eqref{eq:4bis} is $\mathcal{L}^1(\BBF^{n,j})_k=0$ while \eqref{eq:3bis} is $\mathcal{L}^2(\BBF^{n,j})_k=0$.
In order to have more structure, we will require that $\delta_k\BBF$ has the following difference  form:
\begin{equation}\label{fluxF}
	\delta_k\BBF=\widehat{\BBF}_{k+1/2}-\widehat{\BBF}_{k-1/2}\end{equation}
where  $\widehat{\BBF}_{k+1/2}$ depends on $P$ arguments, is consistent with $\BBF$ and uniformly Lipschitz continuous with respect to its arguments. Examples will be given in section \ref{delta}.

Thanks to \eqref{consistant}, we see that
\begin{equation}\label{L2moinsL1}\mathcal{L}^2(\BBF)_k-\mathcal{L}^1(\BBF)_k=\dfrac{\Delta t}{\Delta x} \Lambda A\big ( \delta_k\BBF -\delta_k\BBF^{(0)}\big ),\end{equation}
the important fact is that $\varepsilon$ plays no role here.

We will solve the problem \eqref{eq:3bis} with the following defect correction (DeC) method:
\begin{itemize}
	\item Set, for any $k$,  $\BBF_k^{(0)}=(\BF_k^n, \ldots , \BF_k^n)^T$,
	\item Solve for $p=0, \ldots, M-1$ the problem
	\begin{equation}\label{dec}\mathcal{L}^1(\BBF^{(p+1)})=\mathcal{L}^1(\BBF^{(p)})-\mathcal{L}^2(\BBF^{(p)}),\end{equation}
	\item Set $\BBF^{n+1}=\BBF^{(M)}.$
\end{itemize}
\revB{We remark that the operator $\mathcal L^2$ will never be solved directly as it will be applied to the previously computed iteration $\BBF^{(p)}$. The DeC procedure will converge to the solution of the $\mathcal L^2(\BBF^*)=0$ problem by solving iteratively \eqref{dec}.}
We show that if the problem \eqref{eq:3bis} has a unique solution $\BBF^\star$ and  taking $M=q$, we have a formal error of $\Delta t^q$, i.e.,
for a norm to be defined,
$$\Vert \BBF^{(q)}- \BBF^\star\Vert \leq C\Delta t^q,$$
so that the formal accuracy is the same as solving exactly \eqref{eq:3bis}.
Before doing that, we have first to explain how we solve for $\mathcal{L}^1$ and, hence, \eqref{dec}, then we show the error estimate (and define the proper norm).

\subsubsection{Solution of $\mathcal{L}^1(\BBF)=\mathbf{\mathcal{G}}$ and (\ref{dec}).}
Let us first start with $\mathcal{L}^1(\BBF)=\mathbf{\mathcal{G}}$ for any $\mathcal{G} \in $. Applying $\P$ to this equation, we get, for any $k\in \Z$,
$$
\P\BBF_k=\P\mathcal{G}_k+ \P\BBF_k^{(0)}-\dfrac{\Delta t}{\Delta x} \P C\Lambda \delta_k\BBF^{(0)}= \P\mathcal{G}_k+ \mathbf{\mathcal{K}},$$
with 
\begin{equation*}
	\mathbf{\mathcal{K}}=\P\BBF_k^{(0)}-\dfrac{\Delta t}{\Delta x} \P C\Lambda \delta_k\BBF^{(0)}.
\end{equation*}
\revB{The found equation is explicit for $\P \BBF$ and we can in practice compute this term and use it to obtain the solution of the whole operator.}
Substituting $\P \BBF$ into the Maxwellian, we obtain
$$\BBF_k=\mathcal{G}_k+\BBF_k^{(0)}-\dfrac{\Delta t}{\Delta x} C\Lambda \delta_k\BBF^{(0)}+
\mu\bigg (  A\M\big ( \P\mathcal{G}_k+\mathbf{\mathcal{K}}\big ) -A\BBF_k\bigg )+\mu\ba_0
\otimes (\M(\P\BF^{n,0}_k)-\BF^{n,0}_k) ,$$
\revB{where all the unknown terms $\BBF$ depends only linearly on some coefficients that we can collect on the left hand side}
$$\big (\text{Id}_{q\times q}+\mu A\big ) \BBF_k=\mathcal{G}_k+\BBF_k^{(0)}\!\!\!-\dfrac{\Delta t}{\Delta x}  C\Lambda \delta_k\BBF^{(0)}\!\!\!+
\mu  A\M\big ( \P\mathcal{G}_k+\mathbf{\mathcal{K}}\big )+\mu\ba_0\! \otimes (\M(\P\BF^{n,0}_k)-\BF^{n,0}_k).$$
Now, if $\text{Id}_{q\times q}+\mu A$ is invertible, \revB{we can compute only once and store its inverse to obtain an easy solution of the whole problem, i.e.,}
\begin{subequations}
	\label{solL1}
	\begin{equation}\label{solL1:1}
		\begin{split}
			\BBF_k&=\big (\text{Id}_{q\times q}+\mu A\big )^{-1}\bigg ( \mathcal{G}_k+\BBF_k^{(0)}-\dfrac{\Delta t}{\Delta x}  C\Lambda \delta_k\BBF^{(0)}\bigg )\\
			&\quad +\big (\text{Id}_{q\times q}+\mu A\big )^{-1}\mu A\M\big ( \P\mathcal{G}_k+\mathbf{\mathcal{K}}\big )\\
			&\qquad \quad+
			\big (\text{Id}_{q\times q}+\mu A\big )^{-1}\mu \ba_0 \otimes (\M(\P\BF^{n,0}_k)-\BF^{n,0}_k).
		\end{split}
	\end{equation}
\end{subequations}
\revB{In this way, we are able to find a solution of the system $\mathcal L^1 (\BBF) = \mathcal G $ in a \revA{computationally} explicit way: the source term is split into the linearly implicit part and the Maxwellian evaluated in the previously computed $\P \BBF$. The remaining terms are the explicit right hand side $\mathcal{G}$, the explicit convection term $\delta _k \BBF^{(0)}$ and the explicit part of the high order time integration of the source.} 
\revA{There have been many other works using similar techniques in order to explicitly solve the implicit discretization of the kinetic system \eqref{eq:1:2}, inter alia \cite{coron1991Kinetic,AregbaNatalini,pieraccini2007implicit,Torlo}.} \revB{This DeC approach combining the two operators, has the advantage of being arbitrarily high order without building complicated structures. Indeed, we can use this computationally explicit solution to solve \eqref{dec}.}

Setting 
$$\mathbf{\mathcal G}_k=\mathcal{L}_1(\BBF^{(p)})_k-\mathcal{L}_2(\BBF^{(p)})_k,$$
we can apply \eqref{solL1} directly. However, since the source term discretisation is the same, we have simplifications.
Indeed, after little algebra, \eqref{dec} is rewritten as:
\revD{
	\begin{equation}\label{DeC2}
		\begin{split}
			\BBF^{(p+1)}-\mu A\bigg ( \M(\P\BBF^{(p+1)})-\BBF^{(p+1)}\bigg ) =&\BBF^{(0)}
		-\dfrac{\Delta t}{\Delta x} \Lambda A \delta \BBF^{(p)}-\dfrac{\Delta t}{\Delta x}\ba_0 \otimes \Lambda  \delta \BF^{n,0} +\\
		&\mu\ba_0  \otimes (\M(\P\BF^{n,0})-\BF^{n,0}) 
	\end{split}
	\end{equation}
	and we can apply the same technique. We first apply $\P$,
	\begin{equation}\label{eq:DeCProjection}
		\P\BBF^{(p+1)}=\P\BBF^{(0)}-\dfrac{\Delta t}{\Delta x}\P\Lambda A \delta \BBF^{(p)}-\dfrac{\Delta t}{\Delta x} \ba_0 \otimes \P \Lambda\delta \BF^{n,0},
	\end{equation}
	so we know explicitly $
	\P\BBF^{(p+1)}$, and then we can solve
	\begin{equation*}
	\begin{split}
	\big (\text{Id}_{q\times q} +\mu A\big ) \BBF^{(p+1)}=&\mu A \M(\P\BBF^{(p+1)}) + \BBF^{(0)}
	-\dfrac{\Delta t}{\Delta x} \Lambda A \delta \BBF^{(p)}-\dfrac{\Delta t}{\Delta x}\ba_0 \otimes \Lambda  \delta \BF^{n,0}+\\
	&\mu\ba_0  \otimes (\M(\P\BF^{n,0})-\BF^{n,0}), 	\end{split}
	\end{equation*}
	and then
	\begin{equation}
		\label{solDec}
		\begin{split}
		\BBF^{(p+1)}=\big (\text{Id}_{q\times q} +\mu A\big )^{-1}\bigg (&\mu A \M(\P\BBF^{(p+1)}) + \BBF^{(0)}
		-\dfrac{\Delta t}{\Delta x} \Lambda A \delta \BBF^{(p)}\\
		&-\dfrac{\Delta t}{\Delta x}\ba_0 \otimes \Lambda  \delta \BF^{n,0}+\mu\ba_0  \otimes (\M(\P\BF^{n,0})-\BF^{n,0})\bigg ).
		\end{split}
\end{equation}}
Again, we see that the method is \revA{computationally} explicit.

Next, we show the error estimate, and then we comment more on \eqref{solDec}, in particular when $\varepsilon\rightarrow 0$.
\subsubsection{Error estimate}
If $\varphi: \R\rightarrow \R^p$ is $C^1(\R)^p$ and has a compact support, we can consider the discrete version of its $L^2$ and $H^1$ norms:
$$\Vert \varphi\Vert_{L^2}^2=\sum_{j\in \Z} \Delta x\Vert \varphi_i\Vert^2, \quad \Vert \varphi\Vert_{H^1}^2=\Vert \varphi\Vert_{L^2}^2+\sum_{j\in \Z} \Delta x\Vert D_i\varphi\Vert^2$$ where
$D_i\varphi=\tfrac{\varphi_{i+1}-\varphi_i}{\Delta x}.$

We will establish error estimates that are valid in a given but arbitrary  compact $I=[a, b]$ with discrete equivalent of $L^2_{loc}$  and $H^{-1}_{loc}$ estimates:
$$\Vert \BBF\Vert_{2,I}=\sup\limits_{ \varphi \in C^1_0([a,b])^p}\dfrac{\sum_{j} \Delta x \langle \varphi_i, \BF_i\rangle}{\Vert \varphi\Vert_{L^2}} \text{ and }
\Vert \BBF\Vert_{-1,I}=\sup\limits_{ \varphi \in C^1_0([a,b])^p}\dfrac{\sum_{j} \Delta x \langle \varphi_i, \BF_i\rangle}{\Vert \varphi\Vert_{H^1}} $$
and we note that for $\varphi\in C^1_0([a,b])^p$,  we have a Poincar\'e like inequality
$$\Vert \varphi\Vert _{2,I}\leq (b-a) \Vert D\varphi\Vert _{2,I}.$$

We first show that
\begin{lemma}\label{lemme1}
	If $\widehat{\BBF}_{k+1/2}=\sum\limits_{l=-p}^q \alpha_l \BBF_{k+l}$ and letting $C=\max\limits_{-p\leq l\leq q} \vert \alpha_l\vert \times \max\limits_i\vert \lambda_i\vert$, we have 
	\begin{equation}
		\label{L2-L1}
		\Vert \mathcal{L}^2(\BBF)-\mathcal{L}^1(\BBF)\Vert_{-1,I}\leq C \Vert \BBF\Vert_{2,I} \; \Delta t.
	\end{equation}
\end{lemma}
\begin{proof}
	We have, from \eqref{L2moinsL1} and since $\delta_k\BBF=\widehat{\BBF}_{k+1/2}-\widehat{\BBF}_{k-1/2}$, we have, using that $\varphi$ has a compact support,
	\begin{equation*}
		\begin{split}
			\left \lvert \sum_k \Delta x\langle\varphi_k, \mathcal{L}^2_k(\BBF)-\mathcal{L}^1_k(\BBF)\rangle \right \rvert &=\left \lvert \sum_k \Delta t\langle\varphi_k,A\Lambda\big (\widehat{\BBF}_{k+1/2}-\widehat{\BBF}_{k-1/2}\big )\rangle\right \rvert \\
			&=\left \lvert \sum_k \Delta t \Delta x \langle D_{k+1/2}\varphi, A\Lambda \widehat{\BBF}_{k+1/2}\rangle \right \rvert\\
			&\leq \lVert A \rVert \, \Vert \Lambda\Vert \; \Delta t\sqrt{ \sum_k \Delta x\Vert  D_{k+1/2}\varphi\Vert^2 }\sqrt{ \sum_k \Delta x\Vert\widehat{\BBF}_{k+1/2}\Vert^2}\\
			&\leq C \Delta t\Vert \varphi\Vert_{H^1}  \Vert\BBF\Vert_{2,I}.
		\end{split}
	\end{equation*}
	\revD{We remark that the norm of the coefficients $A$ is smaller or equal to 1.}
\end{proof}
We also have the following lemma on $\mathcal{L}_1$:
\begin{lemma}\label{lemme2}
	We assume that the Maxwellian is Lipschitz continuous and that there exists $C,C'>0$ such that for all $\varepsilon>0$,
	$$\Vert \big (\text{Id}_{(q-1)\times (q-1)}+\mu A\big )^{-1}\Vert \leq C,\quad \mu \Vert \big (\text{Id}_{(q-1)\times (q-1)}+\mu A\big )^{-1}A\Vert \leq C'.$$
	
	Let us consider $\BBF, \BBF'$ and $\mathbf{\mathcal{G}}, \mathbf{\mathcal{G}}'$ such that
	$$\mathcal{L}_1(\BBF)=\mathbf{\mathcal{G}} \text{ and } \mathcal{L}_1(\BBF')=\mathbf{\mathcal{G}}'.$$
	Then, there exists $\alpha>0$, independent of $\BBF, \BBF'$,  $\varepsilon$ and $I$ such that
	$$\Vert \BBF - \BBF'\Vert_{2,I} \leq \alpha \Vert \mathbf{\mathcal{G}}-\mathbf{\mathcal{G}}'\Vert_{2,I}$$
	and
	$$\Vert \BBF - \BBF'\Vert_{-1,I} \leq \alpha \Vert \mathbf{\mathcal{G}}-\mathbf{\mathcal{G}}'\Vert_{-1,I}$$
\end{lemma}
\begin{proof}
	We have the explicit solution $\BBF$ and $\BBF'$ from \eqref{solL1}, and we see that
	$$\big (\text{Id}_{(q-1)\times (q-1)}+\mu A\big )\bigg (\BBF_k-\BBF'_k\bigg )=
	\mathcal{G}_k-\mathcal{G}'_k+\mu A\M\big (\P\mathcal{G}_k-\P\mathcal{G}'_k \big ),$$
	so that if $\Vert ~.~\Vert$ is any of the two norms, we have
	$$\Vert\BBF_k-\BBF'_k\Vert \leq \alpha \Vert \mathcal{G}_k-\mathcal{G}'_k\Vert.$$
	The constant $\alpha$ depends on $C$, $C'$, $\Lambda$ and the Lipschitz constant of the Maxwellian. Remember also that in \eqref{solL1}, $\mathcal{K}$ depends only on $\BBF_0$, $\frac{\Delta t}{\Delta x}$ and $\Lambda$. It is independent of $\varepsilon$.
\end{proof}

Then, wrapping all together,  we have the following proposition:
\begin{proposition}
	Under the assumptions of lemma \ref{lemme1} and \ref{lemme2},  if $\BBF^\star$ is the unique solution of $\mathcal{L}^2(\BBF)=0$,  there exists $\theta$ independent of $\varepsilon$ such that
	we have, for all $p\in \N$
	\begin{equation}\label{inequality}\Vert \BBF^{(p+1)}-\BBF^\star\Vert_{L^2} \leq \big (\theta \Delta t\big )^{p+1} \Vert \BBF^{(0)}-\BBF^\star\Vert_{L^2} .
	\end{equation}
\end{proposition}
\begin{proof}
	We first have, since $\mathcal{L}^2(\BBF^\star)=0$
	\revB{
		\begin{equation*}
			\begin{split}
				\mathcal{L}^1(\BBF^{(p+1)})-\mathcal{L}^1(\BBF^\star)&=\big (\mathcal{L}^1(\BBF^{(p)})-\mathcal{L}^1(\BBF^\star)\big )-\mathcal{L}^2(\BBF^{(p)})\\
				&=\big (\mathcal{L}^1(\BBF^{(p)})-\mathcal{L}^1(\BBF^\star)\big )- \big (\mathcal{L}^2(\BBF^{(p)})-\mathcal{L}^2(\BBF^\star)\big ),
			\end{split}
	\end{equation*}}
	so that combining the inequalities of lemmas \ref{lemme1}, \ref{lemme2} and the Poincar\'e like inequality, we obtain the result.
\end{proof}
\begin{remark}[Comments about inequality \eqref{inequality}]
	This result  shows that  after $p+1$ iteration, the error is $O(\Delta x^{p+1})=O(\Delta t^{p+1}) $ if a CFL-like condition is available. Of course it is better that $\theta \Delta x<1$ for the inequality to be effective, so we may experience a reduction of the CFL number. This reduction needs to be studied  case by case, however this also show that the overall cost of the method is of the order of an explicit one. This is why we name this computationally explicit.
\end{remark}
\revB{\subsubsection{Asymptotic preservation} \label{sec:AP}
	We can show that the presented method is asymptotic preserving (AP), starting from the Chapman--Enskog expansion of the model \eqref{eq:2:1}. Let us define $\bbu^\varepsilon := \P \BF$, we obtain that 
	\begin{equation}\label{eq:chapEnsk}
		\begin{split}
			&\BF = \M(\bbu^{\varepsilon} ) +\mathcal{O}(\varepsilon)  \\
			&\dfrac{\partial \bbu^\varepsilon}{\partial  t} +\frac{\partial \bbf(\bbu^\varepsilon)}{\partial x} = \mathcal{O}(\varepsilon).
		\end{split}
	\end{equation}
	\begin{proposition}
		The discretisation given by \eqref{solDec} is consistent with the limit model \eqref{eq:chapEnsk} up to an $\mathcal{O}(\varepsilon)$.
		\begin{proof}
			Now, using first \eqref{solDec} and then \eqref{eq:DeCProjection}, defining $\bbu^{(p),\varepsilon} = \P \BBF^{(p)}$ and recalling that $\mu = \frac{\Delta t}{\varepsilon}$, $\M (\P \bbu) = \bbu$, $\P \Lambda \M(\bbu) = \bbf(\bbu)$, by induction on the subtimesteps $p$
			$$\bbu^{(p),\varepsilon}=\P \BBF^{(0)} - \frac{\Delta t}{\Delta x}\P \Lambda A \delta \BBF^{(p)}- \frac{\Delta t}{\Delta x}\P \ba_0  \otimes \delta \BBF^{(p)}+\mathcal{O}(\varepsilon),$$
			we can extend the formal expansion also in the discrete case, i.e.,
			\begin{equation}\label{eq:chapEnskDiscrete}
				\begin{split}
					\BBF^{(p+1)} &= \M \left(\P \BBF^{(0)} - \frac{\Delta t}{\Delta x}\P \Lambda A \delta \BBF^{(p)}- \frac{\Delta t}{\Delta x}\P\Lambda \ba_0  \otimes \delta \BBF^{(0)}\right)  +\mathcal{O}(\varepsilon)= \M \left(\bbu^{(p),\varepsilon}\right)  +\mathcal{O}(\varepsilon)\\
					\bbu^{(p+1),\varepsilon} &= \bbu^{(0)} - \frac{\Delta t}{\Delta x} A \delta \P \Lambda \BBF^{(p)} -  \frac{\Delta t}{\Delta x} \ba_0 \otimes\P\Lambda\delta \BBF^{(0)} +\mathcal{O}(\varepsilon)\\
					&=\bbu^{(0)} - \frac{\Delta t}{\Delta x} A \delta \bbf(\bbu^{(p),\varepsilon}) -  \frac{\Delta t}{\Delta x} \ba_0 \otimes\delta\bbf(\bbu^{(0),\varepsilon}) +\mathcal{O}(\varepsilon).
				\end{split}
			\end{equation}
			The final result is a discretisation in space and time of the asymptotic model given by \eqref{eq:chapEnsk}, if the the spatial discretisation is consistent with the space derivative. 
		\end{proof}
		\begin{remark} One can proceed further and prove that, both in the discrete and the continuous case, the next term of the Chapman Enskog expansion is a diffusive term under Whitham's subcharacteristic condition of $\Lambda^2-\partial_u \bbf(u)$ being positive definite. We can also prove that the discretisation is consistent also with that term up to an $\mathcal{O}(\varepsilon^2)+ \mathcal{O}(\Delta t^2)$ if the spatial discretisation is at least consistent. We refer to \cite{Torlo} for the details of such computations for the sake of brevity.
		\end{remark}
	\end{proposition}
	
}

\subsubsection{Examples of $\mathcal{L}^2$ time discretisation.}

Here, we will consider second and fourth order approximation in time in the $\mathcal{L}^2$ operator, namely the Crank-Nicholson method and the fourth order one that uses the points $t_n$, $t_n+\tfrac{\Delta t}{2}$ and $t_{n+1}$. They are described by their matrices $A$,
\begin{itemize}
	\item Second order
	$$A_2=\begin{pmatrix}\frac{1}{2}\end{pmatrix}, \ba_0=\begin{pmatrix}   \frac{1}{2}\end{pmatrix}, \quad \BBF=\begin{pmatrix}\BF^{n,1}\\ \BF^{n,0}\end{pmatrix}$$
	with $\BF^{n,0}=\BF^n$ and $\BF^{n,1}\approx \BF(t^{n+1})$.
	Writing the $\mathcal{L}^2$ discretisation of the time derivative applied to $\BF$, we would have
	$$\BF^{n,1}-\BF^{n,0}+\Delta t A\Lambda \delta\BBF=0,$$ i.e.,
	$$\BF^{n,1}-\BF^{n,0}+\Delta t\big (\frac{1}{2}\Lambda\delta\BF^{n,1} + \frac{1}{2}\Lambda\delta\BF^{n,0}\big)=0.$$
	This is Crank-Nicholson.
	We see that
	$$\big (\text{Id}_{1,1}+\mu A_2\big )^{-1}=\dfrac{2\varepsilon}{2\varepsilon+\Delta t },\qquad  \mu A_2\big (\text{Id}_{1,1}+\mu A_2\big )^{-1}=\dfrac{2\Delta t}{2\varepsilon+\Delta t}.$$
	are uniformly bounded in $[0,2]$ for any $(\Delta t, \varepsilon)$.
	\item The fourth order scheme is obtained by
	$$A_3=\begin{pmatrix}
		\frac{1}{3} & \frac{-1}{24}\\
		\\
		\frac{2}{3} &\frac{1}{6}
	\end{pmatrix},  \ba_0=\begin{pmatrix}\frac{5}{24}  \\\frac{1}{6}\end{pmatrix},  \BBF=\begin{pmatrix}\BF^{n,2}\\ \BF^{n,1} \\ \BF^{n,0}\end{pmatrix}$$
	where $\BF^{n,0}=\BF^n$,   $\BF^{n,1}\approx \BF(t_n+\frac{\Delta t}{2})$ and $\BF^{n,2}\approx\BF(t^{n+1}).$
	We see that
	$$\det\bigg ( \text{Id}_{2\times 2}+\mu A_3\bigg )=\big ( 1+\frac{\mu}{3}\big )\big ( 1+\frac{\mu}{6}\big )+\frac{1}{36}>0$$
	so the matrix is invertible. It is also easy to see that the matrices
	$$\bigg ( \text{Id}_{2\times 2}+\mu A_3\bigg )^{-1} \text{ and }\mu \bigg ( \text{Id}_{2\times 2}+\mu A_3\bigg )^{-1} A$$
	are uniformly bounded.
	
	In fact, in this case, the operator $\mathcal{L}^2$ corresponds to the scheme Lobatto III, which is 4th order accurate \cite{Hairer}. For that reason, we will use this temporal scheme in conjunction with a fourth order spatial approximation.
\end{itemize}
\subsection{Space discretisation: Definition of the $\delta$ operator}\label{delta}

The only question left  is how to define a stable scheme. \revB{As we have seen in section \ref{sec:AP}, the scheme is asymptotic preserving. Under Whitham's subcharacteristic conditions the relaxation terms introduces diffusion which further stabilize the scheme. Hence, we focus only on the stability of the convection scheme, which will also guarantee the stability of the full scheme.} The answer for the fully non linear convection problem is out of reach, at least for this paper, so we will rely on a classical linear stability analysis. The stability of the convection schemes splits into two sub-questions: is the convection scheme defined by $L^2=0$ conditionally or unconditionally stable, and then, is the convection scheme defined by the DeC iteration \eqref{solL1} stable, and under which conditions. In the next section, we will provide 3 examples with increasing accuracy, and sketch a general method.

The matrix $\Lambda$ is diagonal. In \cite{Iserle}, the author considers the transport equation
$$u_t+a u_x=0$$ and shows that if $a<0$ and 
$$u_x(x_i)\approx \frac{1}{\Delta x}\sum_{j=-r}^s \alpha_j u_{i+j},$$
then the order is at most $2\min (r+1,s)$ and in addition the only stable methods are those defined for $r=s$ or $s=r+1$ or $s=r+2$. If $a>0$, we set
$$u_x(x_i)\approx \frac{1}{\Delta x}\sum_{j=-s}^r\alpha_j u_{i+j},$$ while in that case $r=s$ or $r=s+1$ or $r=s+2$. We will only consider these approximations.
Following \cite{Iserle}, we have
\begin{equation*}
	\begin{split}
		\alpha_j&=\dfrac{(-1)^{j+1}}{j}\dfrac{r!s!}{(r+j)!(s-j)!},  \qquad -r\leq j\leq s, j\neq 0,\\
		\alpha_0&=-\sum_{j=-r, j\neq 0}^s \alpha_j
	\end{split}
\end{equation*}
and
\begin{equation*}
	\begin{split}
		\dfrac{\delta_k u}{\Delta x}-\dpar{u}{x}(x_k)&=c\Delta x^q \dfrac{\partial^{q+1}u}{\partial x^{q+1}}(x_k)+O(\Delta x^{q+1}), \qquad q=r+s,\\
		c&=\dfrac{(-1)^{s-1} r!s!}{(r+s+1)!}.
	\end{split}
\end{equation*}
\begin{remark}[Conservation]
	We note that we can always write
	\begin{equation}
		\label{flux}
		\delta_i u=\hf_{i+1/2}-\hf_{i-1/2}
	\end{equation}
	with
	\begin{equation}
		\hf_{i+1/2}=\sum\limits_{j=-r+1}^{s}\beta_j u_{i+j}, \qquad 
		\beta_{j}=\sum_{l\geq j+1}\alpha_l
	\end{equation}
\end{remark}
\begin{proof}
	Assuming that $\hf_{i+1/2}=\sum\limits_{j=-r}^{s-1}\beta_j u_{i+j}$ for any $i$, we write
	\begin{equation*}
		\begin{split}\alpha_{-r}u_{i-r}+\ldots +\alpha_s u_{i+s}=&\bigg ( \beta_{-r+1}u_{i-r+1}+\ldots \beta_{s}u_{i+s}\bigg )-\bigg ( \beta_{-r+1}u_{i-r}+\ldots \beta_{s}u_{i+s-1}\bigg )\\
			=&- \beta_{-r+1}u_{i-r}+ (\beta_{-r+1}-\beta_{-r}) u_{i-1-r}+\ldots \\
			&+(\beta_{l}-\beta_{l-1})u_{i+l-1}+\ldots+\beta_s u_{i+s}
		\end{split}
	\end{equation*}
	so that $\beta_j=-\sum\limits_{l\geq j+1}\alpha_l$, using that $\sum\limits_{l=-s}^r\alpha_l=0$.
\end{proof}

This means that the approximations \eqref{solL1} and \eqref{solDec}, in the limit $\varepsilon\rightarrow 0$ is always conservative since $\Lambda$ is diagonal, and thanks to \eqref{fluxF}.

We list some possible choices for $\delta$:
\begin{itemize}
	\item First order approximation: this is the upwind scheme. If $a>0$, we take $\delta_1 u_j=u_{j}-u_{j-1}$, while if $a<0$, $\delta_1 u_j=u_{j+1}-u_{j}$. If $a=0$, of course $\delta_1 u_j=0$. The flux is
	$$\hf_{j+1/2}=\frac{1}{2}\big ( u_j+u_{j+1}+\text{sign }(a) (u_{j+1}-u_j)\big ), \quad \text{ sign}(a)=\frac{a}{\vert a \vert }.$$
	\item Second order: for $a<0$,
	$$\delta_2 u_j=-\frac{u_{j-1}}{3}-\frac{u_j}{2}+u_{j+1}-\frac{u_{j+2}}{6}$$
	so that
	$$u_x= \dfrac{1}{\Delta x}\bigg (-\frac{u_{j-1}}{3}-\frac{u_j}{2}+u_{j+1}-\frac{u_{j+2}}{6}\bigg )+c\Delta x^3 \dfrac{\partial^4 u}{\partial x^4} +O(\Delta x^4)$$
	with 
	$$c=-\frac{1}{12}.$$ This corresponds to the $[r,r+2]$ approximation with $r=-1$.
	In term of flux, we have (for $a<0$):
	$$f_{j+1/2}=\frac{1}{6}\big ( 2u_{j}+5u_{j+1}-u_{j+2}\big ).$$
	For $a>0$, we have
	$$f_{j+1/2}=\frac{1}{6}\big ( 2u_{j+1}+5u_{j}-u_{j-1}\big ),$$
	so all in all
	$$\hf_{j+1/2}=\frac{1-\text{ sign }(a)}{12}\big ( 2u_{j}+5u_{j+1}-u_{j+2}\big )+\frac{1+\text{ sign }(a)}{12}\big ( 2u_{j+1}+5u_{j}-u_{j-1}\big ).$$
	\item Fourth order:
	if $r=s=2$, and for any $a$
	$$\delta_4^1 u_j=\dfrac{u_{j+2}-u_{j-2}}{12}+2\dfrac{u_{j+1}-u_{j-1}}{3}$$
	hence $$\dpar{u}{x}-\dfrac{\delta_4^1 u}{\Delta x}=c \Delta x^4 \dfrac{\partial^5 u}{\partial x^5}+O(\Delta x^5)$$
	and if $r=1$, $s=3$ and $a<0$, 
	$$\delta_4^2 u=-\dfrac{u_{j-1}}{4}-\dfrac{5}{6}u_j+\dfrac{3}{2}u_{j+1}-\dfrac{u_{j+2}}{2}+\dfrac{u_{j+3}}{12}.$$

	In term of flux, we have: 
	\begin{itemize}
		\item for $\delta_4^1$, 
		$$\hf_{j+1/2}=a\left ( \dfrac{u_{j+2}}{12}+\dfrac{3}{4}u_{j+1}+\dfrac{3}{4}u_j+\dfrac{u_{j-1}}{12}\right)$$
		\item for $\delta_4^2$,
		\begin{equation*}
			\begin{split}
				\hf_{j+1/2}=&\frac{1-\text{ sign }(a)}{2}\left (\dfrac{u_{j+3}}{12}-\dfrac{5}{12}u_{j+2}+\dfrac{13}{12}u_{j+1}+\dfrac{u_j}{4}\right )\\
				+&\frac{1+\text{ sign }(a)}{2}
				\left ( \frac{u_{j+1}}{4}+\frac{13}{12}u_{j}-\frac{5}{12}u_{j-1}+\frac{u_{j-2}}{12}\right ).	
			\end{split}
		\end{equation*}
	\end{itemize}
\end{itemize}

\section{Stability analysis}\label{sec:stab}
We study the stability of the discretisation of the homogeneous problem. Since $\Lambda$ is diagonal, it is enough to look at the scalar conservation problem. We first look at the implicit method \revD{defined by $\mathcal{L}^2 = 0$}, and then at the DeC iteration that is constructed on top of it. This is done by Fourier analysis, we can assume that $a>0$ and the Fourier symbol of $\delta$ is $g$. The table \ref{fouriersymbols} displays the symbols of the operators.
\begin{table}[h]
	\begin{center}
		\footnotesize
		\begin{tabular}{|c||c|}
			\hline
			Operator &Symbol g \\ \hline   &\\
			$\delta_1 $ & $1\!-\!e^{-i\theta}$\\ &\\
			$\delta_2$ & $\dfrac{1}{3}e^{i\theta}\!+\!\dfrac{1}{2}\!-\!e^{-i\theta}\!+\!\dfrac{1}{6}e^{-2i\theta}$  \\ &\\
			$\delta_4^1$ & $i\bigg ( \dfrac{\sin(2\theta)}{6}\!+\!\dfrac{4}{3}\sin \theta\!\bigg )$ \\ &\\
			$\delta_4^2$ & $\dfrac{e^{i\theta}\!}{4}\!+\!\dfrac{5}{6}\!-\!\dfrac{3}{2}e^{-i\theta}\!+\!\dfrac{1}{2}e^{-2i\theta}\!-\!\dfrac{e^{-3i\theta}}{12}$\\ &\\
			\hline
		\end{tabular}
	\end{center}
	\caption{\label{fouriersymbols} List of Fourier symbols.}
\end{table}

The next step is to evaluate the amplification factors of the method, first without DeC iteration, then with DeC iteration.

\subsection{First order in time}
\revD{For a first order scheme the $\mathcal{L}^2$ operator can be written as an implicit Euler method, though being \revA{computationally} explicit, while the DeC iteration, which consists of one step, resemble the explicit Euler method with CFL constrained $0\leq \lambda \leq 1$, where $\lambda = a\Delta t /\Delta x$}. For the $\mathcal L^2=0$ operator, by Fourier transform, we have
$\hu^{n+1}-\hu^n+\lambda g\hu^{n+1}=0$, so that the amplification factor is $G=\dfrac{1}{1+\lambda g}$ which is  of modulus $\leq 1$  if 
$$2\lambda \Re(g)+\lambda^2 \vert g\vert^2\geq 0.$$
If $\lambda\rightarrow 0^+$, we see that $\Re(g)\geq 0$ is a necessary condition, while if $\lambda\rightarrow 0^-$, $\Re(g)\leq 0$. In all cases, $\lambda \Re(g)\geq 0$ is a necessary condition.
Writing $g=a+ib$, and assuming that $\lambda\neq 0$, we see that this condition reads:
$$2 \lambda {a}+{\lambda^2}(a^2+b^2)=(\lambda a+1)^2+\lambda^2b^2-1\geq 0.$$
We also see that
$$(\lambda a+1)^2+\lambda^2b^2\geq (\lambda a+1)^2\geq 1$$
so that $\lambda \Re(g)\geq 0$ is a necessary and sufficient condition for stability.
The Table \ref{tablestabilite} provides the stability condition for the first, second and fourth order schemes. \revD{For the rest of the discussion we consider $a=1$ and in case it is different, one has to rescale $\lambda$ accordingly, as classically done for CFL conditions.}

\subsection{Second order in time}
In that case the $\mathcal{L}^2=0$ scheme reads:
$$u_i^{n+1}-u_i^n+\frac{\lambda}{2} \big ( \delta u_i^n+\delta u_i^{n+1}\big )=0,$$
for which the amplification factor is simply
$$G=\dfrac{ 1-\frac{\lambda}{2} g}{1+\frac{\lambda }{2}g}.$$
We have $|G|\leq 1$ if and only if $$\lambda \Re(g)\geq 0.$$
Again, the Table \ref{tablestabilite} provides the stability condition for the first, second and fourth order schemes.

The DeC iteration is
$$
u_i^{(p+1)}=u_i^n-\frac{\lambda}{2} \big ( \delta u_i^n+\delta u_i^{(p)}\big ),
$$
so that
\begin{equation*}
	\begin{split}
		G_0&=1\\
		G_{p+1}&=1-\frac{\lambda}{2} \big (g +gG_p)
	\end{split}
\end{equation*}
and  we see that
$$G_{p+1}-G=-\frac{\lambda g}{2} \big ( G_p-G)=\bigg ( -\frac{\lambda g}{2}\bigg )^{p+1}\big ( 1-G\big ).$$
\subsection{Fourth order in time}
Here, the $\mathcal L^2=0$  scheme reads: 
\begin{equation*}
	\begin{split}
		u_i^{n+1/2}&-u_i^n+\lambda \bigg ( \frac{5}{24}\delta u_i^n+\frac{1}{3}\delta u_i^{n+1/2}-\frac{1}{24}\delta u_i^{n+1}\bigg )=0\\
		u_i^{n+1}&-u_i^n+\lambda \bigg ( \frac{1}{6}\delta u_i^n+\frac{2}{3}\delta u_i^{n+1/2}+\frac{1}{6}\delta u_i^{n+1}\bigg )=0
	\end{split}
\end{equation*}
so that the Fourier transform gives
$$\begin{pmatrix}\hat{u}^{n+1/2} \\ \hat{u}^{n+1}\end{pmatrix}=G \begin{pmatrix} \hat{u}^n \\ \hat{u}^n \end{pmatrix} $$ with
$$G= 
\begin{pmatrix}
	1+ \frac{\lambda g}{3}  & -\frac{\lambda g }{24}\\
	\frac{2\lambda g}{3}&1+\frac{\lambda g}{6}\end{pmatrix}^{-1} \begin{pmatrix} 1- \frac{5 \lambda}{24}g\\ 1- \frac{\lambda g}{6}\end{pmatrix}=\begin{pmatrix}G_1\\ G_2\end{pmatrix} $$
and we have to look at $\max \{ |G_1|, |G_2|\}\leq 1$ for the calculation of $\hat{u}^{n+1/2}$ and $ \hat{u}^{n+1}$ to be stable.
We have
$$G=\begin{pmatrix} {\frac {-{g}^{2}{\lambda}^{2}+24}{2\,{g}^{2}
			{\lambda}^{2}+12\,\lambda\,g+24}}\\ \noalign{\medskip}{\frac {{g}^{2}{
				\lambda}^{2}-6\,\lambda\,g+12}{{g}^{2}{\lambda}^{2}+6\,\lambda\,g+12}}
\end{pmatrix}.
$$

Then with obvious notations, the DeC iteration is
\begin{equation*}
	\begin{split}
		v_1^{(p+1)}&-u_i^n+\lambda \bigg ( \frac{5}{24}\delta u_i^n+\frac{1}{3}\delta v_1^{(p)}-\frac{1}{24}\delta v_2^{(p)}\bigg )=0\\
		v_2^{(p+1)}&-u_i^n+\lambda \bigg ( \frac{1}{6}\delta u_i^n+\frac{2}{3}\delta v_1^{(p)}+\frac{1}{6}\delta v_2^{(p)}\bigg )=0
	\end{split}
\end{equation*}
The Fourier analysis gives:
\begin{equation*}
	\hv^{(p+1)}=
	\begin{pmatrix}
		1\revD{-}\lambda \theta_0^1 g\\
		1\revD{-}\lambda \theta_0^2 g\end{pmatrix}\hu^n
	\revD{-}
	\lambda g
	\begin{pmatrix}
		\theta_1^1  & \theta_2^1\\
		\theta_1^2& \theta_2^2 \end{pmatrix} \hv^{(p)}, \qquad \begin{pmatrix}\theta_0^1 & \theta_1^1 &\theta_2^1\\\theta_0^1&\theta_1^2 & \theta_2^2\end{pmatrix}=\begin{pmatrix}
		\dfrac{5}{24}& \dfrac{1}{3}&\dfrac{-1}{24}\\ &&\\
		\dfrac{1}{6} & \dfrac{2}{3} &\dfrac{1}{6}\end{pmatrix}
\end{equation*}

The amplification vector, $G_p$ after the $p$-th iteration is defined by
\begin{equation}\label{Jac}
	\begin{split}
		G_0&=\begin{pmatrix}  1\\ 1 \end{pmatrix}\\
		G_{p+1}& =\begin{pmatrix}
			1\revD{-}\lambda \theta_0^1 g\\
			1\revD{-}\lambda \theta_0^2 g\end{pmatrix}\revD{-}\lambda g
		\begin{pmatrix}
			\theta_1^1  & \theta_2^1\\
			\theta_1^2& \theta_2^2 \end{pmatrix} G_p.
	\end{split}
\end{equation}
We note that, setting $\mathbf{\theta}=\begin{pmatrix}
	\theta_1^1  & \theta_2^1\\
	\theta_1^2& \theta_2^2 \end{pmatrix}$, 
$$G_{p+1}-G=(\revD{-}\lambda g)^p  \mathbf{\theta}^p \bigg (\begin{pmatrix}1 \\1 \end{pmatrix}-G\bigg )$$
and $\rho(\mathbf{\theta})=\dfrac{1}{2\sqrt{3}}$. So  using the spectra decomposition of $\theta$ which has two complex and distinct eigenvalues, we have that 
$$\rho(\theta^p)\leq \mu_p=\sqrt {{\dfrac{17}{16}}+\dfrac{\sqrt {241} }{16} } \bigg ( \dfrac{1}{2\sqrt{3}}\bigg )^p.$$
We get finally
$\mu_1=0.4115783562$, $\mu_2=0.1188124373$, $\mu_3=0.03429819635$, $\mu_4=0.009901036444$, hence the convergence is very quick.
\subsection{Summary of the stability analysis}
Combining these expressions with the actual form of the Fourier symbol of $\delta$, we get the results of table \ref{tablestabilite}.
\begin{table}[h]
	\begin{center}
		\begin{tabular}{|c||ccc|}
			\hline
			& First order & second order & Fourth order\\
			\hline
			$\delta_1$ &      $\checkmark$            &         $\checkmark$                &       $\checkmark$           \\
			$\delta_2$ &           $\checkmark$    &        $\checkmark$               &   $\checkmark$ if $\lambda\leq 4.5$\\ 
			$\delta_4^1$ &        $\checkmark$      &         $\checkmark$ ($|G|=1$)               & $\checkmark$ ($|G|=1$)\\
			$\delta_4^2$ &        $\checkmark$        &          $\checkmark$             & $\checkmark$ if $\lambda\leq \frac{9}{4}$\\
			\hline
			Analytical condition &$\lambda \Re(g)\geq 0$ & $\lambda \Re ( g)\geq 0$ & $\lambda \Re \big ( g-\frac{\lambda g^2}{6}\big ) >0$\\
			\hline
		\end{tabular}
	\end{center}
	\caption{\label{tablestabilite} Stability conditions for the original scheme.}
\end{table}

Now, we turn our attention on the DeC iteration. For the second order in time approximation, we first have
$$G_p=(1-\theta_p)G+\theta_p \qquad \text{with } \theta_p=\big (-1\big )^p\bigg (\dfrac{\lambda}{2}g\bigg)^p.$$
so we get
$$|G_p|^2-1= |1-\theta_p|^2 \big ( |G|^2-1\big ) +2\Re\bigg (  \overline{\theta_p}(1-\theta_p)(G-1)\bigg )$$
hence if $|G|\leq 1$,  a sufficient condition is that
$$\Re\bigg (  \overline{\theta_p}(1-\theta_p)(G-1)\bigg )\leq 0.$$

For the fourth order scheme, we have similarly 
$$\revD{G_p=\big ( \text{Id}-(-\lambda g\theta)^p\big ) G+ (-\lambda g\theta)^p e}, \qquad e=\begin{pmatrix} 1 \\1 \end{pmatrix},$$but it is more complicated to get an analytical condition. So we rely on Maple.

The stability conditions are summarised in Table \ref{stabilityDeC}.
\begin{table}[h]
	\begin{center}
		\begin{tabular}{|cc||c|c|c|c|c|c|}
			\hline
			\multicolumn{2}{|c|}{Scheme} & \multicolumn{6}{c|}{$\#$ iterations}\\
			\hline
			Order & $\delta$ & 1 & 2 & 3 & 4 & 5 &6\\
			\hline
			2 & $\delta_1$& 1 &1&1&1&1&1\\
			2&$\delta_2$ & 0 & $\geq 0.85$& $\geq 1.22$ & $\geq 1.02$& $\geq 1.08$& $\geq 1.23$\\
			2& $\delta_4^1$ & 0 & 0 & $\geq 1.45$ & $\geq 1.45$  & $\geq 0.002$ & $\geq 0.01$\\
			2& $\delta_4^2$ & 0 &$\geq 0.5$ &$\geq 0.69$ & $0.71$ & $0.73$ & $0.73$ \\
			\hline \hline
			3 & $\delta_1$ & $6$ & $\geq 1.5$ & $\geq 1.87$ & $\geq 2$ & $\geq 2.23$ & $\geq 2.48$\\
			3 & $\delta_2$ & 0& 0& $1$ & $\geq 2.0447$ &$\geq 2.17120$ & $\geq 2.568$\\
			3& $\delta_4^1$ & 0 & 0 & 0 & $\geq 1.6171$ & $\geq 2.4727$ & $\geq 2.9162$\\
			3& $\delta_4^2$ & 0 & 0& $\geq 0.1$& $\geq 1.3096$ & $\geq 1.3955$ & $\geq 1.8282$\\
			\hline
		\end{tabular}
	\end{center}
	\caption{\label{stabilityDeC}CFL number for stability of the DeC iterations. $0$ means that the scheme is unconditionally unstable. If a  real number $x$ is given, it  means that the scheme is stable up to CFL $x$, if $\geq x$ is written, this  means that the scheme is stable for at least CFL $x$ (and slightly above indeed).}
\end{table}

%
\section{Wave model}
We have to specify the diagonal matrix $\Lambda$ and the Maxwellians $\M$.
We will use two kinds of wave models:
\begin{itemize}
	\item A two waves model. In that case,
	$$\Lambda=\begin{pmatrix}a&0\\0&-a\end{pmatrix}$$
	with $a\geq \max_i\rho(\bbf'(u_i))$.
	Setting $\bu^\varepsilon=\P\BF$, we have $\bbf(\bu^\varepsilon) = \P\Lambda \BF$ and we know explicitly $\M=(\M_1,\M_2)$:
	\begin{equation}
		\label{Maxwell2}
		\M_1(\P\BF)=\frac{1}{2}\bigg ( \bbu^\varepsilon+\frac{\bbf}{a}\bigg ), \qquad \M_2(\P\BF)=\frac{1}{2}\bigg ( \bbu^\varepsilon-\frac{\bbf}{a}\bigg ).
	\end{equation}
	\item A three waves model, where
	$$\Lambda=\begin{pmatrix}a & 0 & 0\\0&0&0\\0&0&-a\end{pmatrix}.$$
	In that case, the Maxwellian is $\M=(\M_1,\M_2,\M_3)$ and we have
	\begin{equation*}
		\begin{split}
			\bu^\varepsilon& =\;\,\M_1+\M_2+\M_3\\
			\bbf(\bu^\varepsilon)&=a\M_1\qquad-a\M_3
		\end{split}
	\end{equation*}
	so we need to specify $\M_2$. 
\end{itemize}
For the scalar problems, we will use the two wave model that reveals itself sufficient. For the fluid problems, we will show that the two wave models is not perfect, and hence the three wave model needs to be considered.

In the case of 3 waves, let us specify $\M_2$. Following \cite{Bouchut}, we know that the sub-characteristic condition is equivalent to the monotonicity of the Maxwellians: they need to be differentiable and have only positive eigenvalues. In \cite{Bouchut,Natalini}, it is proposed to use 
\begin{equation}
	\begin{split}
		\M_1(\bu^\varepsilon)&= \frac{1}{a}\bbf_+(\bu^\varepsilon)\\
		\M_2(\bu^\varepsilon)&=\bu-\dfrac{\bbf_+(\bu^\varepsilon)-\bbf_-(\bu^\varepsilon)}{a}\\
		\M_3(\bu^\varepsilon)&=\frac{1}{a}\bbf_-(\bu^\varepsilon)\\
	\end{split}
\end{equation}
where $\bbf(\bbu^\varepsilon)=\bbf_+(\bu^\varepsilon)+\bbf_-(\bu^\varepsilon)$, $\bbf_{\pm}$ are differentiable, $\nabla_\bu\bbf_+(\bu)$ has only positive eigenvalues, while $\nabla_\bu\bbf_-(\bu^\varepsilon)$ has only negative eigenvalues.
A possible choice, inspired by the Enquist-Osher-Solomon flux, is
$$\M_2(\bu^\varepsilon)=\dfrac{\int_0^{\bu^\varepsilon} \vert \bbf'(s)\vert \; ds}{\vert a\vert },$$but  the integral (or the path integral for system) must be evaluated.
In the case of the Euler equations, we give a second one  that  does not necessitate to evaluate an integral.
%
In the case of the Euler equations, 
\begin{equation}\label{eq:euler}
	\bbu^\varepsilon=\begin{pmatrix}
		\rho\\ \rho u,\\E\end{pmatrix}, \quad \bbf(\bbu^\varepsilon)=\begin{pmatrix}
		\rho u\\ \rho u^2+p\\ u(E+p)\end{pmatrix}, \quad p=(\gamma-1) \big ( E-\frac{1}{2}\rho u^2\big ),
\end{equation}
we propose to use a  Maxwellian that relies on the van Leer flux splitting \cite{vanLeer}. It is purely algebraic and defined by:
\begin{enumerate}
	\item if $M=\dfrac{u}{c}\leq -1$, with $c^2=\gamma\frac{p}{\rho}$, then $\bbf_-(\bbu^\varepsilon)=\bbf(\bbu^\varepsilon)$, $\bbf_+(\bbu^\varepsilon)=0$,
	\item If $M\geq 1$, then $\bbf_+(\bbu^\varepsilon)=\bbf(\bbu^\varepsilon)$, $\bbf_-(\bbu^\varepsilon)=0$,
	\item if $-1\leq M\leq 1$, then
	$$
	\bbf_-(\bbu^\varepsilon)=\begin{pmatrix}
		Q\\
		\frac{QR}{\gamma}\\
		\frac{QR^2}{2(\gamma^2-1)}
	\end{pmatrix}, \; Q=-\dfrac{\rho}{4c}(u-c)^2, \; R=(\gamma-1)u-2c
	,$$ and $\bbf_+=\bbf-\bbf_-$
\end{enumerate}
The eigenvalues of $\bbf_{\pm}$ are bounded by
$$a=\left \{\begin{array}{cc}
	(|u|+c)\dfrac{\gamma+3}{2\gamma+|M|(3-\gamma)} & \text{ if } |M|\leq 1\\
	|u|+c & \text{ else.}
\end{array}\right . $$
Note that $\tfrac{\gamma+3}{2\gamma+|M|(3-\gamma)}\leq \tfrac{\gamma+3}{2\gamma}$ for $|M|\leq 1$. For $\gamma=1.4$, $\tfrac{\gamma+3}{2\gamma}=\tfrac{11}{7}\approx 1.57$.

\section{Non linear stabilisation}\label{stabilisation}
If the solution is expected to be non smooth, then one can expect the occurrence of spurious oscillations. Sometimes, oscillations are acceptable, provided they do not lead to the crash of the simulation. In order to get rid of them, or to control them, we have adopted the MOOD technique initially designed in \cite{Clain} with some improvements described in \cite{Vilar}. We have adapted it our way in order to get results that are formally of order $p+1$ in space and time, here $p=1,2,3$.

\revB{MOOD is an \textit{a posteriori }corrector of high order numerical methods \cite{Clain,Vilar}. MOOD requires a sequence of schemes ordered from the most accurate/less stable one to the low order/more reliable one. It also requires a series of criteria that the solution should fulfill, e.g. physical admissibility, discrete minimum principle or numerical errors. After having performed a step of the most accurate scheme, it checks the criteria on each cell/degree of freedom, and detects the areas where the criteria are not met. 
	There, we switch to the next scheme in a \textit{cascade} style, which is supposed to be more stable and reliable. We proceed iteratively until either the criteria are met or the most reliable/less accurate \textit{parachute} scheme is used. The \textit{parachute} scheme should analytically guarantee all the criteria.
	
	In the following we describe how the criteria must be verified on the described spatial discretisation, while the list of the schemes that we use consists always of 2 schemes (the considered one and the upwind discretisation $\delta_1$ as \textit{parachute} scheme) and we specify directly in the numerical simulations which criteria will be considered, as they are problem dependent.}

We proceed as follows: at  the time step $t_n$, we have the values $(\BF_k^{n})_{k}$. For now on, we drop the superscript $n$, since there is no ambiguity.
In the DeC iteration \eqref{eq:3}, with the spatial scheme defined by $\delta_p$, writes (with the convention that $\BF^{(l),0}=\BF^0$ for $l=0, \ldots , q-1$)
\revD{
	\begin{equation*}
		\BF_k^{(p+1),j}-\BF_k^0+\dfrac{\Delta t }{\Delta x}\bigg (\sum_{l=0}^q a_{jl}\Lambda \delta_k \BF^{(p),l}-\mu \sum_{l=0}^q a_{jl}\big ( \M\P\BF_k^{(p+1),l}-\BF_k^{(p+1),l}\big )\bigg ), \; p=0, \ldots , q-1,
	\end{equation*}
	from which we get
	\begin{equation}\label{update}
		\P\BF_k^{(p+1),j}-\P\BF_k^0+\dfrac{\Delta t }{\Delta x}\bigg (\sum_{l=0}^q a_{jl}\P\Lambda \delta_k \BF^{(p),l}\bigg )=0, \quad p=0, \ldots , q-1.
	\end{equation}
	The increment $\delta_k \BF^{l}$ is the difference of two terms, and we write 
	$$\delta_k \BF^{l}=\Lambda\big ( \hat{\BF}_{k+1/2}^l-\hat{\BF}_{k-1/2}^l\big )= \Phi_{k}^{[k, k+1],l}+\Phi_{k}^{[k-1,k],l}$$
	with
	$$\Phi_{k}^{[k, k+1],l}= \Lambda \hat{\BF}_{k+1/2}^l-\Lambda \BF_{k}^l, \quad \Phi_{k+1}^{[k,k+1],l}=\Lambda \BF_{k+1}^l-\Lambda\hat{\BF}_{k+1/2}^l.$$
	One equivalent way to rephrase the conservation is 
	\begin{equation}
		\label{conservation}
		\Phi_{k}^{[k, k+1],l}+\Phi_{k+1}^{[k,k+1],l}=\Lambda\big ( \BF_{k+1}^l-\BF_{k}^l\big )
	\end{equation}
	and the right hand side of this relation is independent of the order $p$.  It is equivalent because we see that
	$$\Lambda \hat{\BF}_{k+1/2}^l=\frac{1}{2}\bigg ( \Lambda \BF_{k+1}^l+\Lambda \BF_{k}^l-\big ( \Phi_{k}^{[k, k+1],l}-\Phi_{k+1}^{[k,k+1],l}\big ) \bigg ).$$
	Using this we rewrite \eqref{update} as:
	\begin{subequations}\label{update:rd}
		\begin{equation}
			\label{update:rd:1}
			\P\BF_k^{(p+1),j}= \frac{1}{2}\bigg ( \big (\widetilde{\P\BF_k^{(p+1),j}}\big )_{k-1/2}+ \big (\widetilde{\P\BF_k^{(p+1),j}}\big )_{k+1/2}\bigg )
		\end{equation}
		with
		\begin{equation}
			\label{update:rd:2}
			\begin{split}
				\big (\widetilde{\P\BF_k^{(p+1),j}}\big )_{k-1/2}=\P\BF_k^{0}-\dfrac{\Delta t }{\Delta x}  \Phi_{k}^{[k-1,k],(p),j}\\
				\big (\widetilde{\P\BF_k^{(p+1),j}}\big )_{k+1/2}=\P\BF_k^0-\dfrac{\Delta t }{\Delta x}  \Phi_{k}^{[k,k+1],(p),j}
			\end{split}
		\end{equation}
	\end{subequations}
	In practice, we compute for each interval $[k, k+1]$
	\begin{equation}
		\label{update:rd:3}
		\begin{split}
			\big (\widetilde{\P\BF_k^{(p+1),j}}\big )_{k+1/2}=\P\BF_k^0-\dfrac{\Delta t }{\Delta x}  \Phi_{k}^{[k,k+1],(p),j}\\
			\big (\widetilde{\P\BF_{k+1}^{(p+1),j}}\big )_{k+1/2}=\P\BF_{k+1}^0-\dfrac{\Delta t }{\Delta x}  \Phi_{k+1}^{[k,k+1],(p),j}
		\end{split}
	\end{equation}
	and then apply \eqref{update:rd:1}.}
%

In the simplified version of the MOOD algorithm we use, we  consider only two spatial approximations, namely the first order one defined by $\delta_1$, and the high order one define by $\delta_p$, $p=2$ or $3$ in this paper. The idea is to use as often as possible the highest order scheme, and to use the low order one  to correct potential problems. Knowing the $\{\BF_k^{(p),j}\}_k$, we first compute  for each interval $[k,k+1]$ the quantities defined by \eqref{update:rd:3} with the high order residuals. Then we test the values of the results using a set of criteria, applied on $(\widetilde{\P\BF_k^{(p+1),j}})_{k+1/2}$ and $(\widetilde{\P\BF_{k+1}^{(p+1),j}})_{k+1/2}$. This set of criteria  is explained in the next paragraph. If both $(\widetilde{\P\BF_k^{(p+1),j}})_{k+1/2}$ and $(\widetilde{\P\BF_{k+1}^{(p+1),j}})_{k+1/2}$ pass the tests, this element is declared sane, else un-sane. This enable to identify a set $\mathcal{I}$ of  un-sane elements $[k,k+1]$ where the criteria are not met, and we store the residual $\{ \Phi_{k}^{[k,k+1],(p),j}, \Phi_{k+1}^{[k,k+1],(p),j}\}$ for the sane elements.  
We then repeat the procedure for the un-sane elements with the lowest order scheme. 
At the end of the procedure, we have evaluated  residuals, that we still denote by $\{ \Phi_{k}^{[k,k+1],(p),j}, \Phi_{k+1}^{[k,k+1],(p),j}\}$, even though they are potentially evaluated by different schemes.
We then compute $\P\BF_k^{(p+1),j}$ by \eqref{update}. There is no problem of conservation since \eqref{conservation} holds true.

\bigskip
Now we describe the criteria we apply to $(\widetilde{\P\BF_k^{(p+1),j}})_{k+1/2}$ and $(\widetilde{\P\BF_{k+1}^{(p+1),j}})_{k+1/2}$, following the ideas of \cite{Clain,Vilar} with some small adaptation to the context.  When specific tests are done on a variable, we denote this variable by $\xi$. For a scalar problem, $\xi$ is simply the conserved variable. In the case of the Euler equations, we test this on some primitive variables: the density and the energy, and for some severe problems, the velocity. We can add as many criteria as needed.
\begin{enumerate}
	\item 
	We first check  if $(\widetilde{\P\BF_k^{(p+1),j}})_{k+1/2}$ and $(\widetilde{\P\BF_{k+1}^{(p+1),j}})_{k+1/2}$ lies in the invariance domain if relevant: in the case of the Euler equation, we check if the density and the internal energy are both positive. If not, we set the criteria to .FALSE. on this element. In that case we jump to the next element, else we look for the next criterion.
	\item We check if the solution is not locally constant. Taking $\nu=\Delta x^3$ and $\mathcal{S}$ the stencil defined by the operator $\delta$, we check if 
	$$\big |\max_{l\in \mathcal{S}}\xi_{i+l}-\min_{l\in \mathcal{S}}\xi_{i+l}\big |\leq \nu \text{ and } \big |\max_{l\in \mathcal{S}}\xi_{i+1+l}-\min_{l\in \mathcal{S}}\xi_{i+1+l}\big |\leq \nu$$
	If this is true, the criteria is kept to .TRUE., else it is set to .FALSE. and we jump to the next element
	\item We check if  a new extrema is created or not, by comparing with the solution at the previous time step, in a neighbourhood extended to the right and the left by one cell: we are running at CFL 1. 
	\revD{
		\begin{enumerate}
			\item We first test if $\xi_k^{n+1},\xi_{k+1}^{n+1}\in [\min\limits_{l\in \mathcal{S}}\xi_{k+l}+\epsilon, \max\limits_{l\in \mathcal{S}}\xi_{k+l} -\epsilon]\cap  [\min\limits_{l\in \mathcal{S}}\xi_{k+1+l}+\epsilon, \max\limits_{l\in \mathcal{S}}\xi_{k+1+l} -\epsilon]$.
			If this is true, we jump to the next element, 
			\item else, denoting by $P_j$ the Lagrange interpolation polynomial that interpolates $\{\xi_{j+l}\}_{l\in \mathcal{S}}$ 
			\begin{itemize}
				\item we compute $\xi'=P_k'(x_k)$, $\xi'_L=P'_k(x_k-\frac{\Delta x}{2})$, $\xi'^{,k-1/2}_{min/max}=\min/\max \big (P'_k(x_k-\frac{\Delta x}{2}), P'_{k-1}(x_k-\frac{\Delta x}{2})\big )$ then
				\begin{itemize}
					\item If $\xi'_L<\xi'$, $\alpha_L=\min(1, \dfrac{\xi'^{,k-1/2}_{max}-\xi'}{\xi'_L-\xi'})$
					\item if $\xi'_L=\xi'$, $\alpha_L=1$
					\item if $\xi'_L<\xi'$, $\alpha_L=\min(1, \dfrac{\xi'^{,k-1/2}_{min}-\xi'}{\xi'_L-\xi'})$
				\end{itemize}
				\item We compute $\xi'=P_{k}'(x_{k})$, $\xi'_R=P'_k(x_k+\frac{\Delta x}{2})$, $\xi'^{,k+1/2}_{min/max}=\min/\max \big (P_{k+1}'(x_k+\frac{\Delta x}{2}), P'_{k}(x_k+\frac{\Delta x}{2})\big )$ then
				\begin{itemize}
					\item If $\xi'_R<\xi'$, $\alpha_R=\min(1, \dfrac{\xi'^{,k+1/2}_{max}-\xi'}{\xi'_R-\xi'})$
					\item if $\xi'_R=\xi'$, $\alpha_R=1$
					\item if $\xi'_R<\xi'$, $\alpha_R=\min(1, \dfrac{\xi'^{,k+1/2}_{min}-\xi'}{\xi'_R-\xi'})$
				\end{itemize}
				\item We set $\alpha=\min(\alpha_L,\alpha_R)$
				\item If $\alpha=1$, then we have a true extrema, keep the criteria to .TRUE. and jump to the next element. Else, we set the criteria to .FALSE. and jump to the next element.
			\end{itemize}
	\end{enumerate}}
	
\end{enumerate}
\revD{
	The idea behind the step 3  is described in \cite{Vilar} and is also related to \cite{kuzmin}: we try to check if the gradient of the interpolation $\xi$ lies in the interval $\left[\min \big (\xi'^{,k-1/2}_{min},\xi'^{,k+1/2}_{min}\big ), \max\big (\xi'^{,k-1/2}_{max},\xi'^{,k+1/2}_{max}\big )\right]$.
}

\revB{
	\begin{remark}[Stability]
		The von Neumann stability study of Section \ref{sec:stab} does not hold directly on the MOOD algorithm but it is clear that we are dealing with a combination of the stabilities of the schemes used in the MOOD \textit{cascade}. Hence, if we choose CFL conditions that guarantee the stability of both the high order schemes and of the \textit{parachute} scheme (upwind CFL=1 in our case), then we know that the global MOOD scheme will be von Neumann stable.
\end{remark}}
\section{Numerical examples}
\subsection{Scalar problems}
The first problem is the transport equation with periodic boundary conditions
\begin{equation*}
	\dpar{u}{t}+\dpar{u}{x}=0
\end{equation*}
where the initial condition is 
\begin{equation}\label{IC}
	u_0(x)=\sin(2\pi x)+0.5.
\end{equation} A two waves model is used with $a=1.01$ (so a little larger that the actual maximum speed. We always proceed as such for scalar and system cases.
We make a convergence test for short and long final times, namely $T=0.5$ and $T=10$.  The CFL number, with respect to the wave model maximum speed, 
is always set to $1$. In both cases, we see that  the expected order of accuracy is obtained. Note that \revD{for the second and fourth order schemes, the non linear stabilisation does not detect any troubled point, giving exactly the same error as in the non stabilised case.} We show the results for the fourth order schemes. All the calculations are done with the two waves model. 
Note that the first order scheme, with the initial condition given by the Maxwellian, is nothing more than the Lax-Friedrichs scheme, for second order in time approximation. For fourth  order in time, since the equilibrium relaxation is more complex, we get a different scheme.  Note that the non linear stabilisation procedure of section \ref{stabilisation} does not flag any cell.
\begin{table}[h]
	\begin{center}
		\begin{tabular}{||c||cccccc|}
			\multicolumn{7}{c}{First order}\\
			\hline
			$h$ & $L^1$ & $r$ & $L^2$ & $r$ & $L^\infty$ & $r$\\
			\hline
			50&      $ 2.75963992\;10^{-2}$&-& $ 3.89822088\;10^{-2}$&-&   $2.49972343\;10^{-2}$&-\\
			100&    $  1.31966826\;10^{-2 }$&1.43& $ 1.86553914\;10^{-2 }$&1.43&  $1.18893785\;10^{-2}$&1.43\\
			200&   $  6.94037229\;10^{-3 }$&1.33& $ 9.81052034\;10^{-3}$&1.33&   $6.26028096\;10^{-3}$&1.33\\
			400&   $    3.47535103\;10^{-3}$&1.38&   $4.91373939\;10^{-3}$&1.38&  $ 3.13191721\;10^{-3}$&1.38\\
			800&   $  1.73895701\;10^{-3  }$&1.38& $2.45900149\;10^{-3}$&1.38&   $1.56634545\;10^{-3}$&1.38\\
			\hline
			\multicolumn{7}{c}{Second order}\\ 
			%
			\hline
			50   &   $ 4.83627617\;10^{-3  }$&-& $6.70434069\;10^{-3  }$&-& $4.40502120\;10^{-3}$&-\\
			100    &  $ 1.21754361\;10^{-3} $ &2.07&$ 1.70489808\;10^{-3 }$ &2.06& $1.10206485\;10^{-3}$&2.08\\
			200   &  $   3.05118738\;10^{-4 } $&2.07&$ 4.29379230\;10^{-4}$  &2.07&$ 2.75470491\;10^{-4}$&2.08\\
			400    & $   7.60697367\;10^{-5 } $&2.08& $1.07314205\;10^{-4 }$ &2.08&$ 6.85840860\;10^{-5}$&2.08\\
			800    &   $ 1.89899602\;10^{-5}$  &2.08& $2.68216700\;10^{-5 } $&2.08& $1.71091069\;10^{-5}$&2.08\\
			\hline
			\multicolumn{7}{c}{Fourth order}\\
			\hline
			
			50  &    $ 0.201424\;10^{-4  } $&-& $0.278979\;10^{-4 } $ &-& $0.183213\;10^{-4}$&-\\
			100  &      $0.122376\;10^{-5  }$& 4.04& $0.171337\;10^{-5}$  &4.02& $ 0.110818\;10^{-5}$&4.04\\
			200   &   $0.758547\;10^{-7 }$  &4.01&$0.106742\;10^{-6}$ &4.00&$0.684850\;10^{-7}$&4.01\\
			400   &   $ 0.472475\;10^{-8 } $ &4.00& $0.666515\;10^{-8}  $&4.00&  $0.425979\;10^{-8}$&4.00\\
			800   &    $ 0.294831\;10^{-9} $  &4.00& $0.416433\;10^{-9} $ &4.00& $ 0.265631\;10^{-9}$&4.00\\
			\hline
			\multicolumn{7}{c}{Fourth order+MOOD}\\
			\hline
			50  &    $ 0.201424\;10^{-4  } $&-& $0.278979\;10^{-4 } $ &-& $0.183213\;10^{-4}$&-\\
			100  &      $0.122376\;10^{-5  }$& 4.04& $0.171337\;10^{-5}$  &4.02& $ 0.110818\;10^{-5}$&4.04\\
			200   &   $0.758547\;10^{-7 }$  &4.01&$0.106742\;10^{-6}$ &4.00&$0.684850\;10^{-7}$&4.01\\
			400   &   $ 0.472475\;10^{-8 } $ &4.00& $0.666515\;10^{-8}  $&4.00&  $0.425979\;10^{-8}$&4.00\\
			800   &    $ 0.294831\;10^{-9} $  &4.00& $0.416433\;10^{-9} $ &4.00& $ 0.265631\;10^{-9}$&4.00\\
			\hline
		\end{tabular}
		\caption{\protect\label{T=0.5} Order of convergence for the convection problem and two wave model for order 1, 2 and 4, and 4th order with MOOD. The final time is $T=0.5$. One can see that the two fourth order results are \emph{identical} as expected.}
	\end{center}
\end{table}
\begin{table}[h]
	\begin{center}
		\begin{tabular}{||c||cccccc|}
			\multicolumn{7}{c}{First order}\\
			\hline
			$h$ & $L^1$ & $r$ & $L^2$ & $r$ & $L^\infty$ & $r$\\
			\hline
			50 &       $0.374521524 $ &-&   $ 0.529551387 $ &-&   $ 0.337824076 $   &-\\
			100  &      $0.222015068 $ &1.22&   $ 0.313846916$  &1.22&    $0.200003594 $   &1.21\\
			200 &     $  0.121430904 $  &1.23&  $ 0.171709701 $   &1.30& $ 0.109344706 $   &1.3\\
			400  &     $  6.35740533\;10^{-2 } $&1.34&$ 8.99051651\;10^{-2} $&1.34& $ 5.72395548\;10^{-2}$&1.34\\
			800  &      $ 3.25359367\;10^{-2} $ &1.36& $4.60121371\;10^{-2 }$&1.36&  $2.92929020\;10^{-2}$&1.36\\
			\hline
			\multicolumn{7}{c}{Second order}\\ 
			\hline   

			50 &    $    9.76886451\;10^{-2}$ &-& $ 0.135240585 $   &-&   $ 8.88576061\;10^{-2}$&-\\
			100  &     $  2.43498404\;10^{-2 }$&2.08&  $ 3.40924263\;10^{-2}  $&2.07&  $2.20487341\;10^{-2}$&2.09\\
			200  &     $  6.06694631\;10^{-3 } $&2.08& $ 8.53730459\;10^{-3 }$ &2.08&  $5.47759095\;10^{-3}$&2.08\\
			400  &      $ 1.51354610\;10^{-3}  $&2.08&  $2.13514664\;10^{-3 }$&2.08&  $ 1.36459176\;10^{-3}$&2.08\\
			800  &     $  3.77953198\;10^{-4 }$ &2.08& $ 5.33837010\;10^{-4 }$&2.08&   $3.40518804\;10^{-4}$&2.08\\
			\hline
			%
			%
			%
			\hline
			\multicolumn{7}{c}{Fourth order}\\
			\hline

			50    &    $ 0.399627\;10^{-3}$ &-&  $ 0.554329\;10^{-3 }$ &-& $ 0.363964\;10^{-3}$&-\\
			100  &   $  0.244527\;10^{-4} $&4.03& $ 0.342394\;10^{-4 } $&4.01& $ 0.221427\;10^{-4}$&4.03\\
			200  &     $  0.151613\;10^{-5}$ &4.01&  $ 0.213344\;10^{-5} $&4.00&   $0.136893\;10^{-5}$&4.01\\
			400  &     $  0.944521\;10^{-7}$ &4.00&  $0.133241 \;10^{-6 }$&4.00&  $ 0.851587\;10^{-7}$&4.00\\
			800   &    $  0.589189\;10^{-9}$ &4.00&  $ 0.832214\;10^{-8 }$&4.00& $  0.530836\;10^{-8}$&4.00\\
			\hline
			\multicolumn{7}{c}{Fourth order with Mood}\\
			\hline

			
			50    &    $ 0.399627\;10^{-3}$ &-&  $ 0.554329\;10^{-3 }$ &-& $ 0.363964\;10^{-3}$&-\\
			100  &   $  0.244527\;10^{-4} $&4.03& $ 0.342394\;10^{-4 } $&4.01& $ 0.221427\;10^{-4}$&4.03\\
			200  &     $  0.151613\;10^{-5}$ &4.01&  $ 0.213344\;10^{-5} $&4.00&   $0.136893\;10^{-5}$&4.01\\
			400  &     $  0.944521\;10^{-7}$ &4.00&  $0.133241 \;10^{-6 }$&4.00&  $ 0.851587\;10^{-7}$&4.00\\
			800   &    $  0.589189\;10^{-9}$ &4.00&  $ 0.832214\;10^{-8 }$&4.00& $  0.530836\;10^{-8}$&4.00\\
			\hline
		\end{tabular}
		\caption{\protect\label{T=10} Order of convergence for the convection problem and two wave model for order 1, 2 and 4 with MOOD. The final time is $T=10$. The fourth order results with and without stabilisation are \emph{identical} as expected.}
	\end{center}
\end{table}

The convergence tables are obtained against the solution of the asymptotic model, proving moreover that the model is asymptotic preserving as expected.

\begin{figure}[h]
	\begin{center}
		\subfigure[First order, second and fourth order solution, with no stabilisation.]{\includegraphics[width=0.45\textwidth]{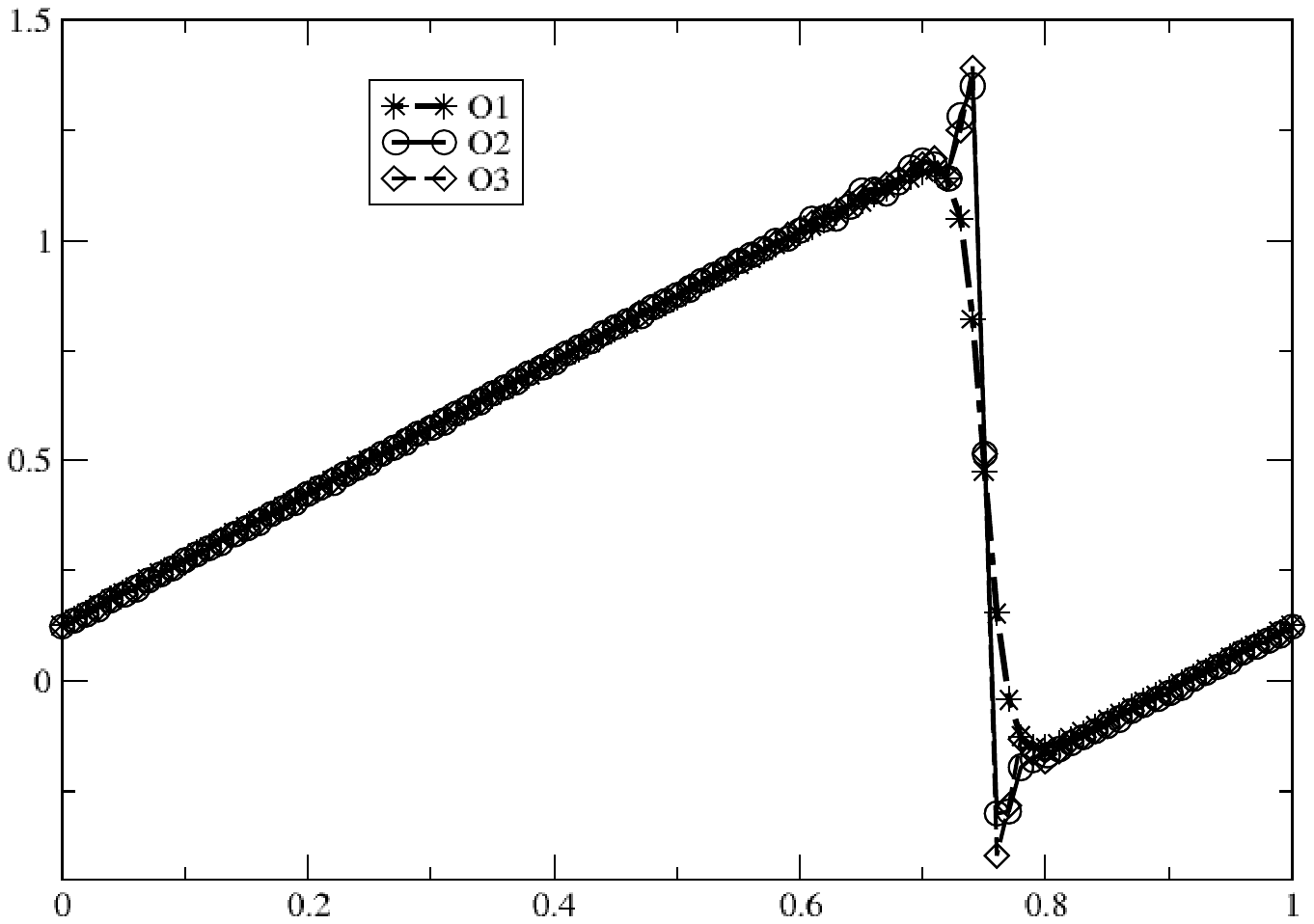}} \hfill \subfigure[First order, second and fourth order with MOOD.]{\includegraphics[width=0.45\textwidth]{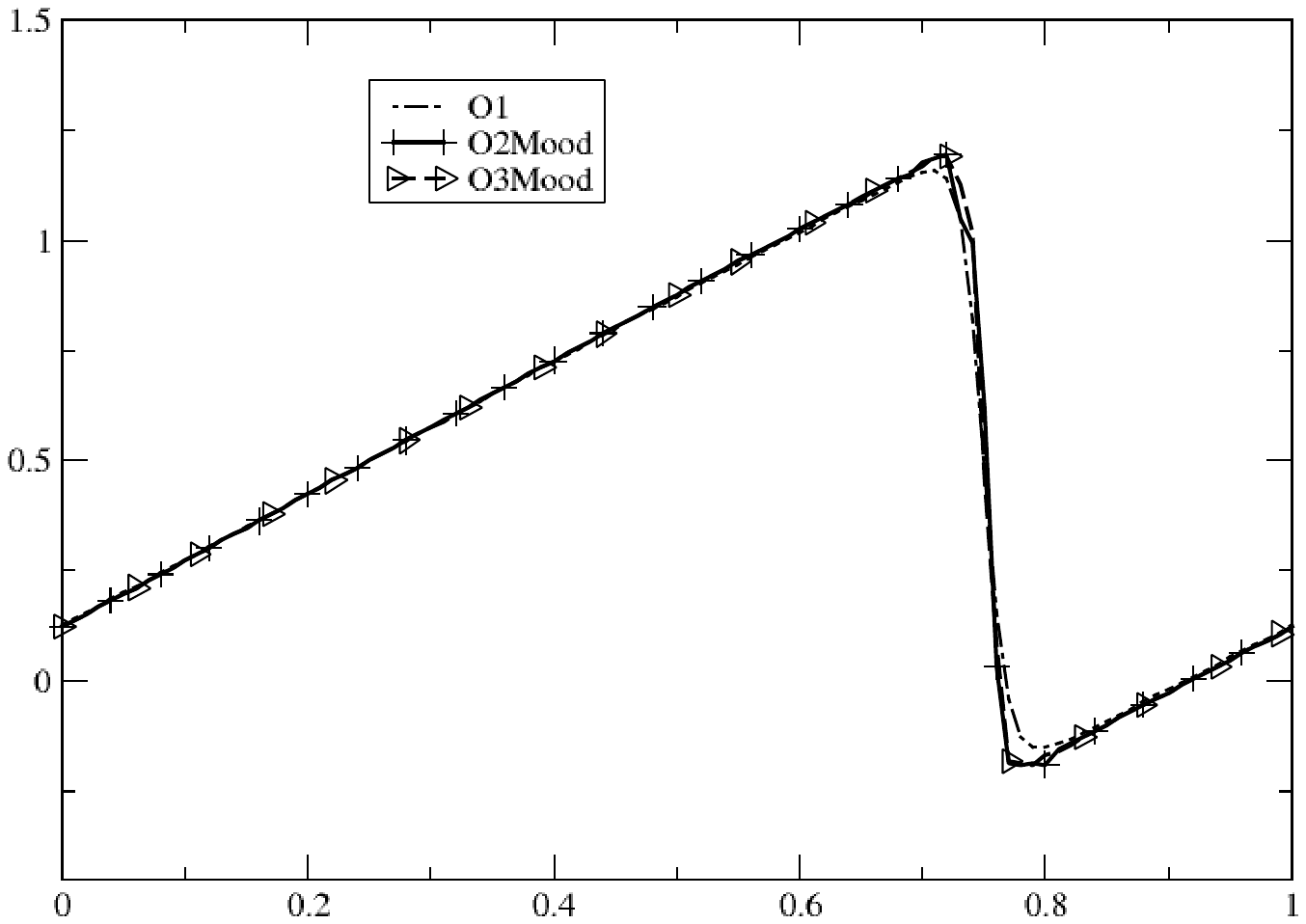}}
	\end{center}
	\caption{Burgers equation T=0.5. Initial condition : $u_0(x)=\sin( 2 \pi x)+0.5$}\label{burger}
\end{figure}

The figure \ref{burger} shows some results for the Burgers equation
\begin{equation}\label{burg}
	\dpar{u}{t}+\frac{1}{2}\dpar{u^2}{x}=0
\end{equation}
with the initial condition \eqref{IC}.  This generates an unsteady shock wave, so a priori more challenging than a steady one. The non linear stabilisation performs correctly.

The last scalar example is the Buckley-Leverett equation
\begin{equation}\label{buck}
	\dpar{u}{t}+\dpar{f(u)}{x}=0, \qquad f(u)=\dfrac{u^2}{u^2+(1-u)^2}
\end{equation}
again with the same initial condition \eqref{IC}. The flux is non convex, so the problem is a bit more challenging.

\begin{figure}[h]
	\begin{center}
		\subfigure[unlimited O2-O4 solutions]{\includegraphics[width=0.45\textwidth]{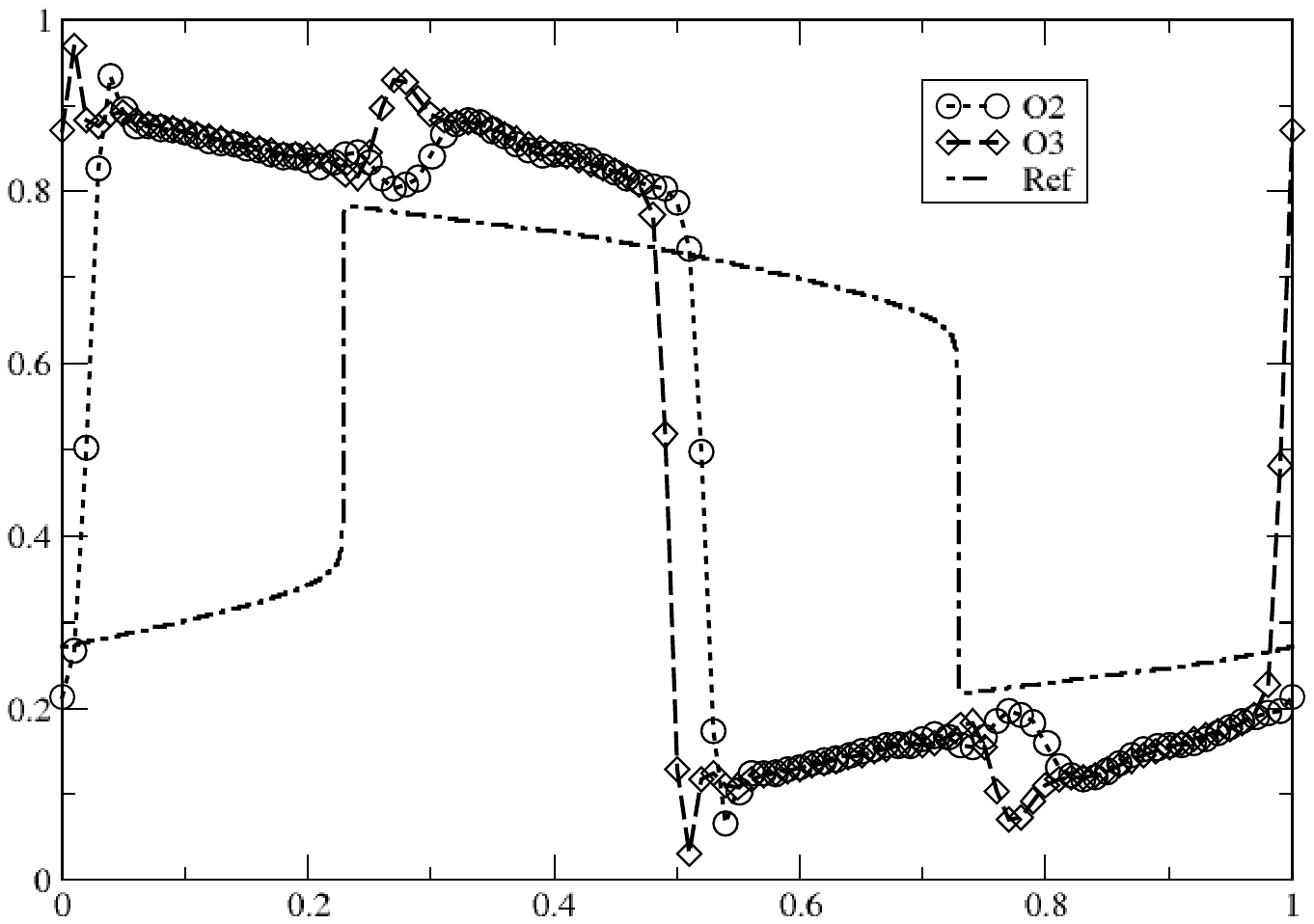}}
		\subfigure[first order and O2-O4 limited solutions]{\includegraphics[width=0.45\textwidth]{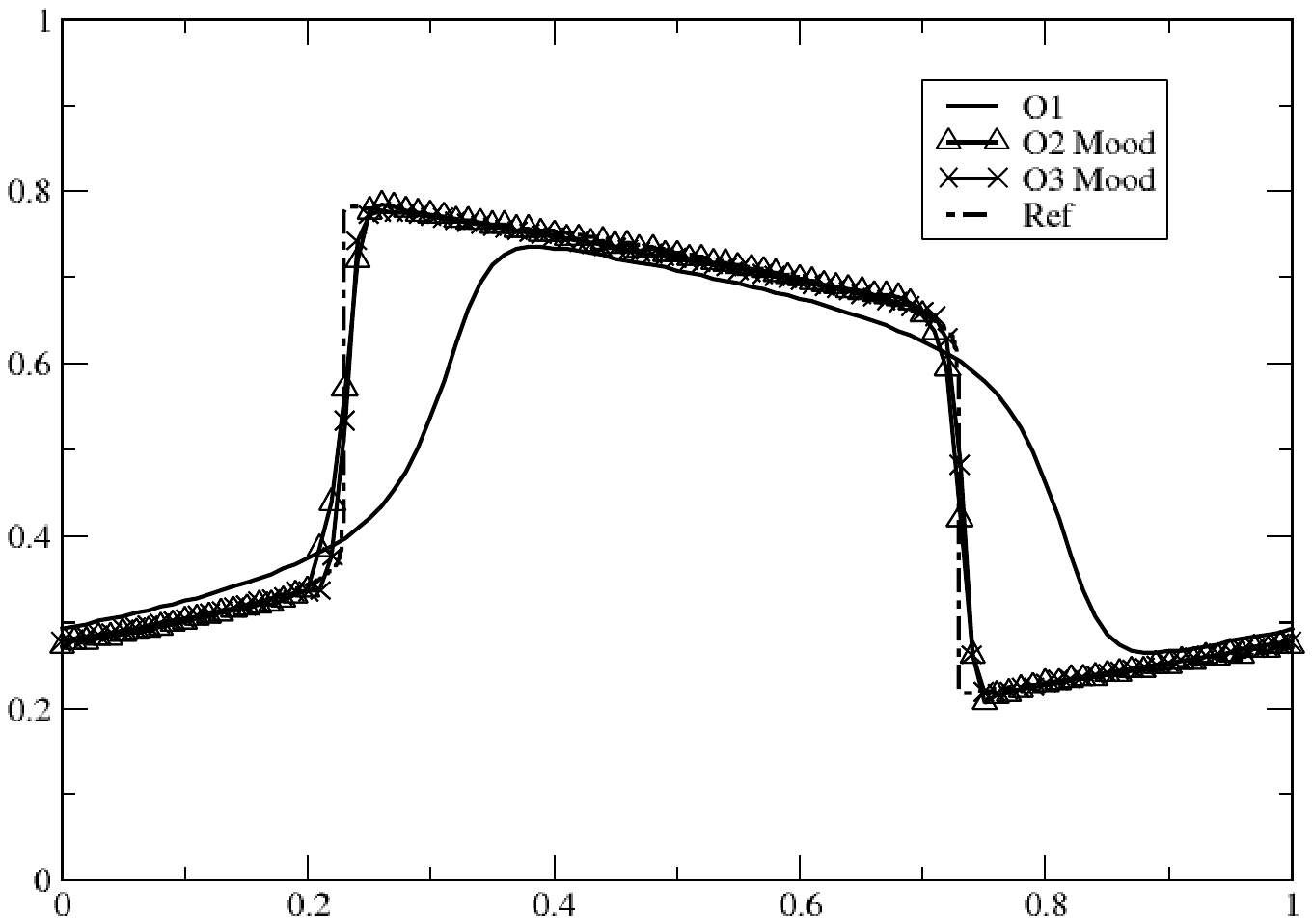}}
	\end{center}
	\caption{ Buckley-Leverett problem with 100 points.  A reference solution (first order with 10000 points) is also indicated.}\label{buck:1}
\end{figure}

We have run the simulations with 100 spatial points, until time $T=1$, with the 2 waves model we have considered. The first order (O1), second order (O2), fourth order (O4), second order with non linear stabilisation (O2M), and fourth order with non linear stabilisation (O4M) are displayed in figure \ref{buck:1}, together with a reference solution computed with 1000 points and the first order scheme: remember that this corresponds to the Lax Friedrichs scheme, and it satisfies all entropy inequalities. This guaranties that the scheme converges. The non linearly stabilized solution have a correct behavior.

\subsection{Uniform order of accuracy with respect to $\varepsilon$}
\begin{table}[h]
	\begin{center}\footnotesize
		\begin{tabular}{|c||c|c||c|c||c|c||c|c||c|c|}
			\hline
			$\varepsilon$ &\multicolumn{2}{c||}{0}  & \multicolumn{2}{c||}{$10^{-6} $}  &      \multicolumn{2}{c||}{$10^{-4} $} & \multicolumn{2}{c||}{$10^{-3}$}& \multicolumn{2}{c|}{$10^{-2}$}\\
			\hline \hline
			$\log \Delta x$ & $\log L^2$ & slope& $\log L^2$ & slope&  $\log L^2$ & slope& $\log L^2$ & slope&$\log L^2$ & slope\\
			\hline
			-2.995   &              -1.332      & - &       -1.332   &        -          &  -1.332       &          -         & -1.400     &-& -2.049       & - \\
			-3.688   &            -2.655    &1.908 &         -2.655    &  1.908               &  -2.655      &   1.907                  &  -2.710  &1.889&-3.314       &  1.825  \\
			-4.382   &            -4.030     &1.984 &        -4.031     &   1.984          &   -4.031      &      1.984             &  -4.063  &1.951& -4.647           & 1.922  \\
			-5.075  &           -5.415     & 1.997   &     -5.415        &  1.997  & -5.415     &     1.996            &  -5.410 &1.943& -5.999         & 1.951  \\
			-5.768   &            -6.801    &  1.999  &      -6.801       &  1.999      &   -6.800                 & 1.998   &-6.749&1.930&-7.365                  & 1.971  \\
			\hline
		\end{tabular}
	\end{center}
	\caption{\label{AP3}Convection problem: error for the second order scheme with different $\varepsilon$}
\end{table}

\begin{table}[h]\begin{center}
		\footnotesize
		\begin{tabular}{|c||c|c||c|c||c|c||c|c||c|c|}
			\hline
			$\varepsilon$ &\multicolumn{2}{c||}{0}  & \multicolumn{2}{c||}{$10^{-6} $}  &      \multicolumn{2}{c||}{$10^{-4} $} & \multicolumn{2}{c||}{$10^{-3}$}& \multicolumn{2}{c|}{$10^{-2}$}\\
			\hline \hline
			$\log \Delta x$ & $\log L^2$ & slope& $\log L^2$ & slope&  $\log L^2$ & slope& $\log L^2$ & slope&$\log L^2$ & slope\\
			\hline
			-2.995   &     -3.615        & -  & -3.615       &   - &  -3.701   &          -         & -3.701    &-  &   -4.445    & - \\
			-3.688   &     -6.385          &3.995   &  -6.385       & 3.996   &  -6.475    &      3.996            &  -6.475   &  4.001  &  -7.212      & 3.992   \\
			-4.382   &     -9.158         & 4.000  &    -9.159      & 4.001   &  -9.253     &     4.000              &   -9.253  & 4.007  &   -9.991    & 4.008 \\
			-5.075  &       -11.93         & 4.000  &    -11.93      &  3.999  &  -12.03    & 4.000                   & -12.03    &  4.007 &   -12.76     & 4.004  \\
			-5.768   &      -14.70         & 4.000  &    -14.70      & 4.000  &  -14.80     &  4.000                 &  -14.81    &  4.006  & -15.53        & 4.001  \\
			\hline
	\end{tabular}\end{center}
	\caption{\label{AP4}Convection problem: error for the fourth order scheme with different $\varepsilon$}
\end{table}

\revD{In order to test the convergence of the scheme also in non asymptotic regimes, we use again the scalar advection equation}
\begin{equation}
	\label{convection}
	\dpar{u}{t}+\dpar{u}{x}=0
\end{equation}
with initial condition 
\begin{equation}\label{sine}
	u(x,0)=\sin (2\pi x).
\end{equation}
The CFL number is set to $1$, the time and space order are set to 2 and 4, the final time is $T=1$, the boundary conditions are periodic. The values of $\varepsilon$ are $\{0, 10^{-6}, 
10^{-5},   10^{-4}, 10^{-3},  10^{-2}\}$. We evaluate the order of convergence with the following standard procedure: if $u_{\Delta x}$, $u_{\Delta x/2}$  are the numerical solution evaluated for consecutive meshes, the order $\alpha$ is, for the norm $\Vert ~.~\Vert$,
$$\alpha=\dfrac{\log \Vert u_{\Delta x}-u_{\Delta x/2}\Vert}{\log \Delta x}.$$
\revD{This test does not have obvious result, as order reduction phenomena are common for IMEX schemes when the space discretization and the relaxation variable are of the same order. Nevertheless, we see on tables  \ref{AP3} and \ref{AP4} that the convergence order does not depend on $\varepsilon$.} It is optimal. 
Also qualitatively, we see that for small enough values of $\varepsilon$, the solution obtained for $\varepsilon=0$ is almost indistinguishable from $0<\varepsilon\ll 1$.

\subsection{Euler equations}
\revD{In this section we test our scheme on Euler equations \eqref{eq:euler}}. We set $\gamma=1.4$ and we run some standard cases: the Sod case and the Shu-Osher case.
\subsubsection{Sod test case}
\begin{figure}[h]
	\begin{center}
		\subfigure[Density]{\includegraphics[width=0.45\textwidth]{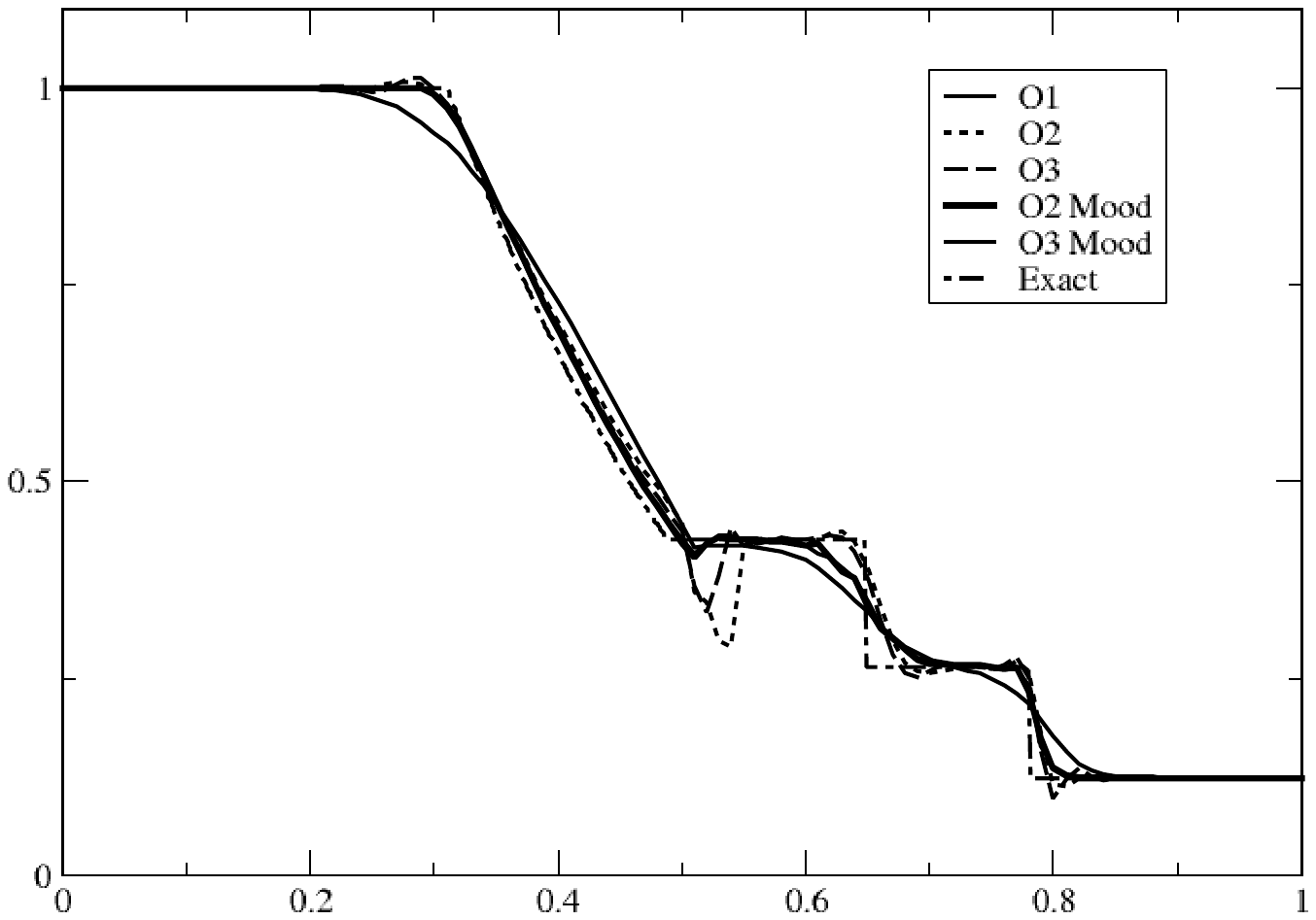}}
		\subfigure[Pressure]{\includegraphics[width=0.45\textwidth]{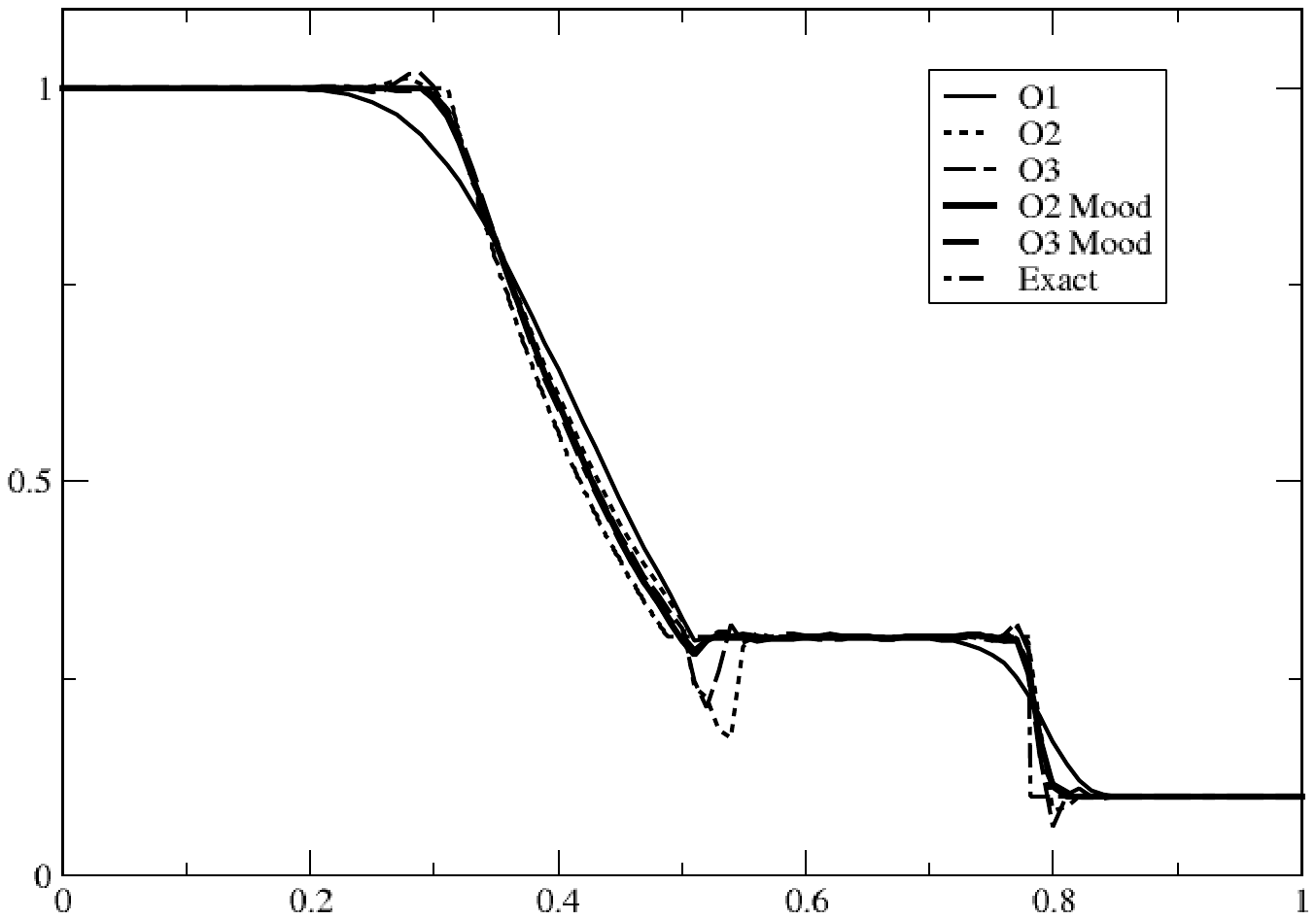}}
		\subfigure[Velocity]{\includegraphics[width=0.45\textwidth]{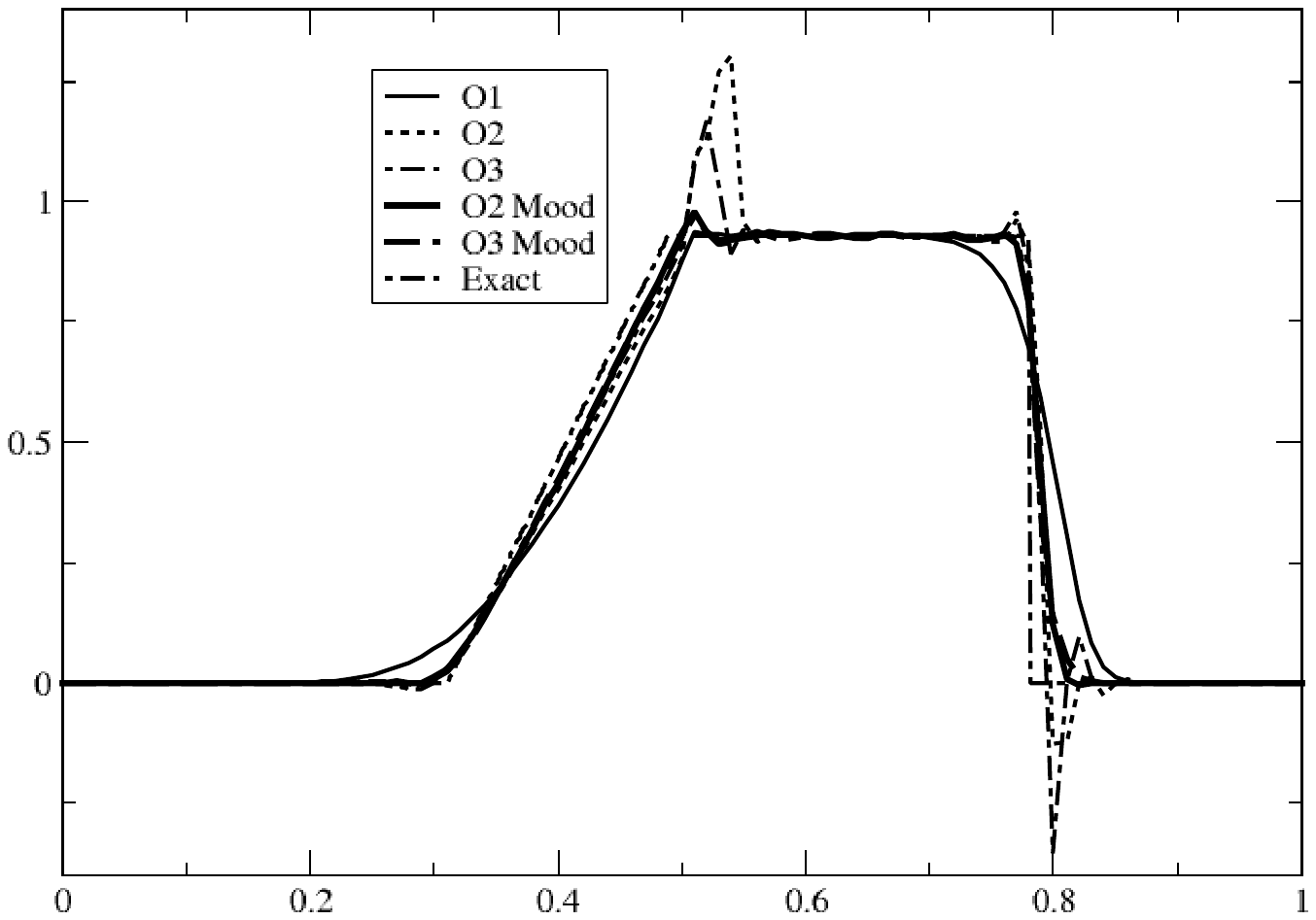}}
	\end{center}
	\caption{ Sod problem with 3-waves model: plot of the density, velocity and the pressure. Displayed solutions for order 1, 2 and 4 schemes with 100 points with and without MOOD. Also the exact solution is plotted}\label{Sod:3w}
\end{figure}

The Sod problem consists of a Riemann problem defined by the following initial conditions:
$$(\rho,u,p)^T=\left \{ \begin{array}{ll}
	(1,0,1)^T & \text{ for } x<0.5m\\
	(0.125,0,0.1)^T &\text{ else.}
\end{array}\right .
$$
The final time is $T=0.16$.
We have used the 3-waves model described above. The mesh resolution is of 100 elements, and the CFL is again 1 in all cases. From figure \ref{Sod:3w}, we see that the results are of good quality, at least compared with more standard methods.

\begin{figure}[h]
	\begin{center}
		\subfigure[Density]{\includegraphics[width=0.45\textwidth]{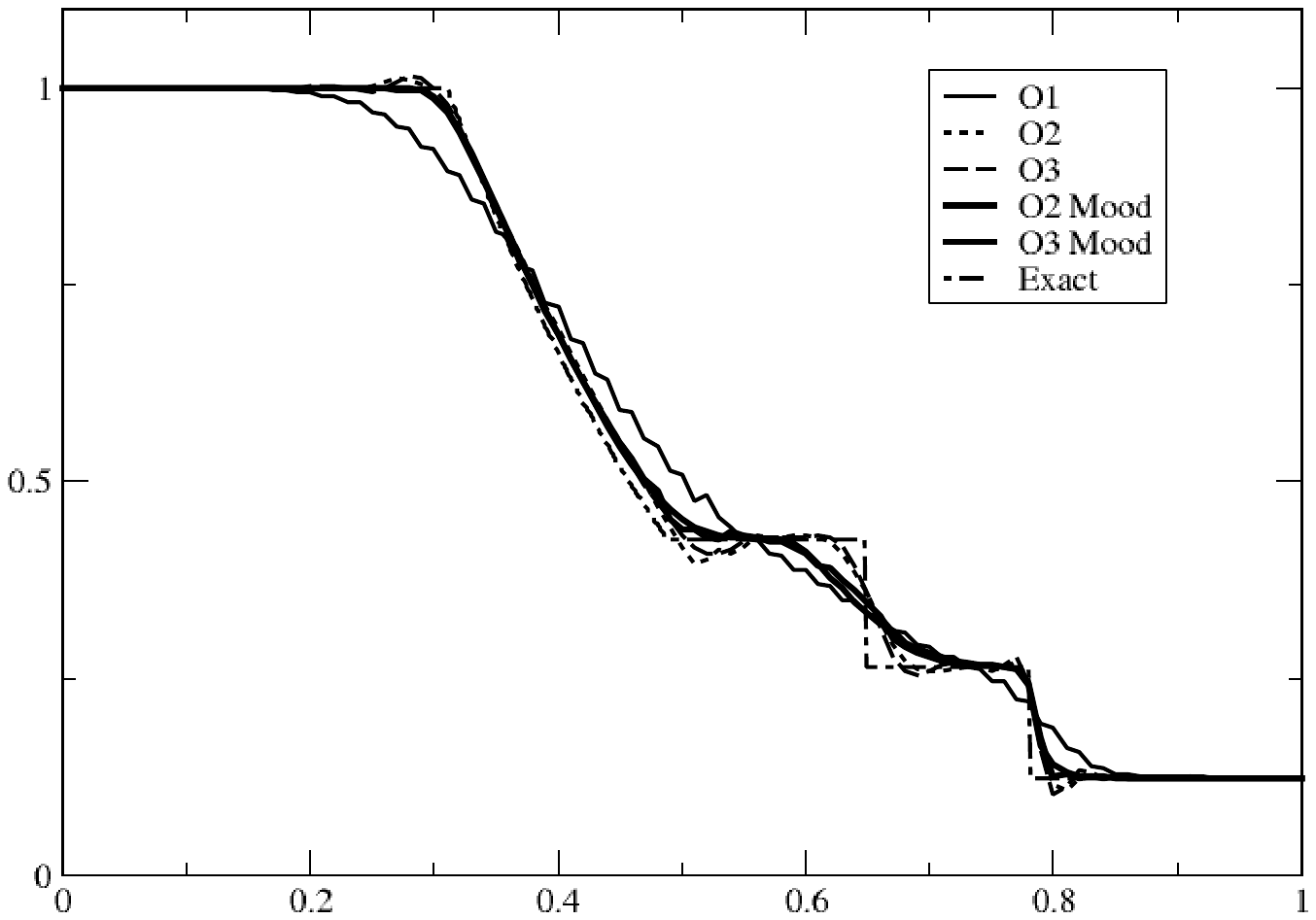}}
		\subfigure[Pressure]{\includegraphics[width=0.45\textwidth]{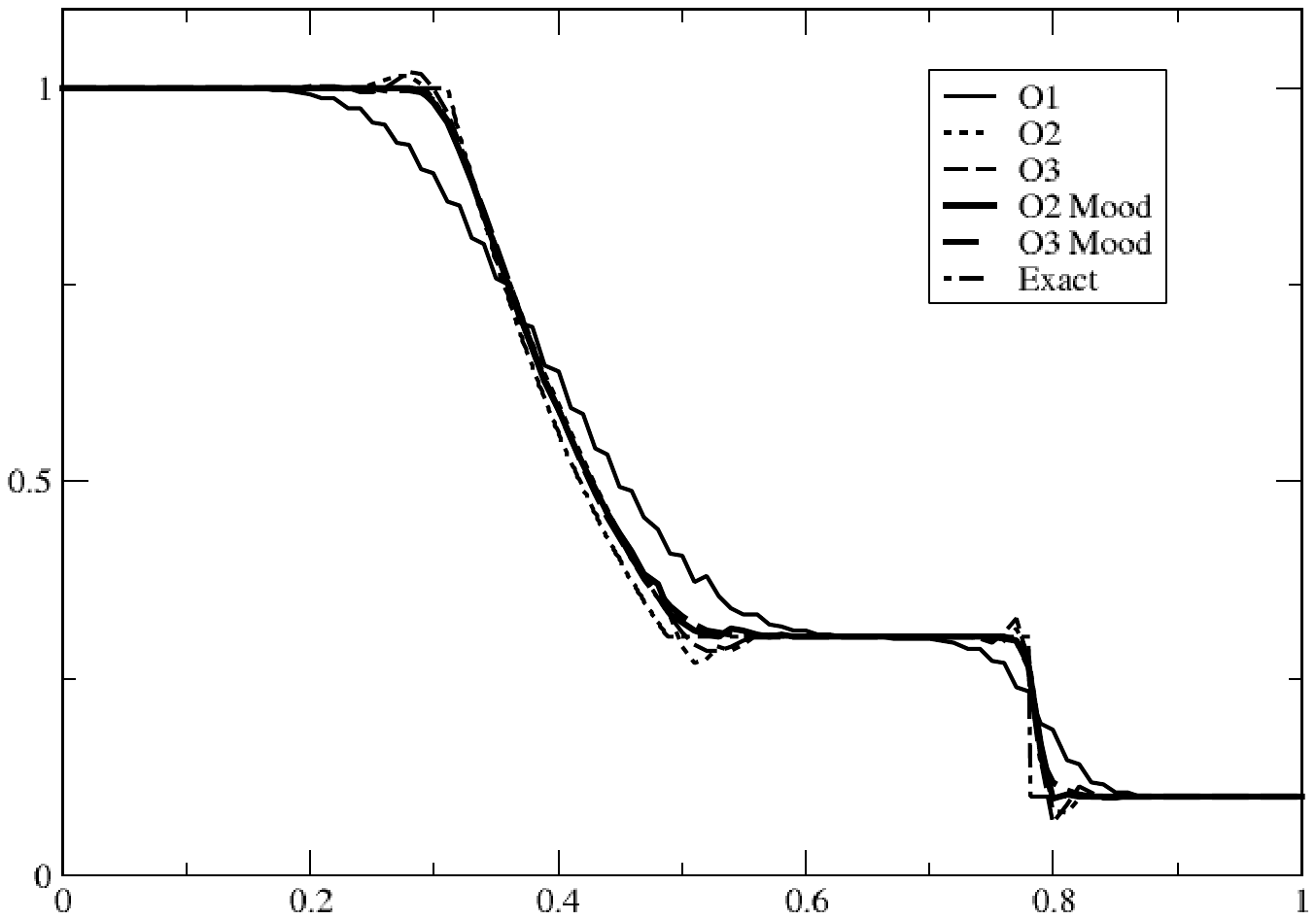}}
		\subfigure[Velocity]{\includegraphics[width=0.45\textwidth]{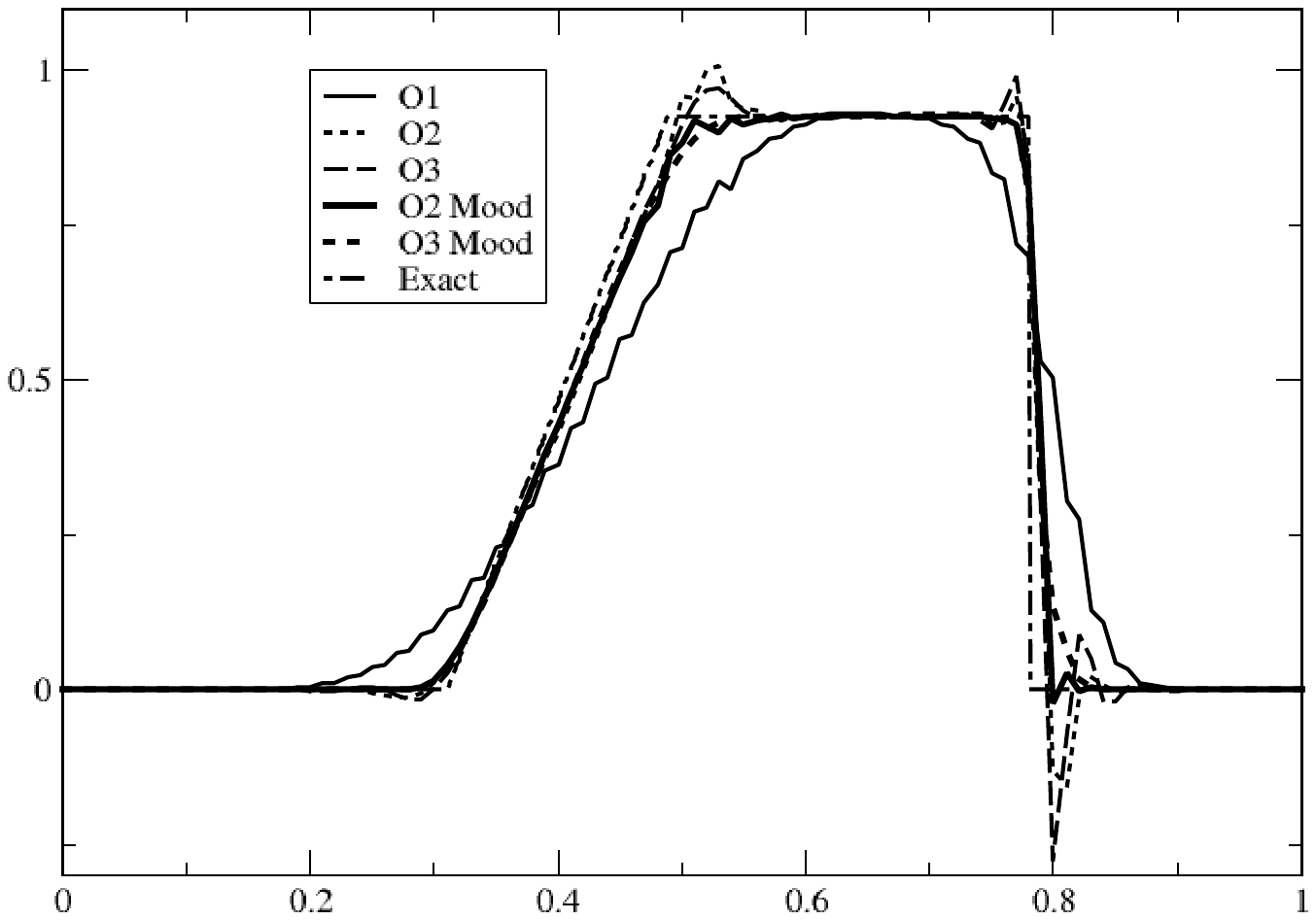}}
	\end{center}
	\caption{ Sod problem with 2-waves model: plot of the density, velocity and the pressure.Displayed solutions for order 1, 2 and 4 schemes with 100 points with and without MOOD. Also the exact solution is plotted}\label{Sod:2w}
\end{figure}

For the sake of completeness, we have made the same simulation with the two wave model in figure \ref{Sod:2w}.
We see a stair case solution of the first order in space which is typical for the Lax-Friedrichs scheme. Comparing the solutions, the 2 waves model provides results of lower quality with respect to the 3 waves one. For that reason, we will not consider it anymore for the Euler equations.

\subsubsection{Shu-Osher test case}
\begin{figure}[h]
	\begin{center}
		\subfigure[O2Mood]{\includegraphics[width=0.45\textwidth]{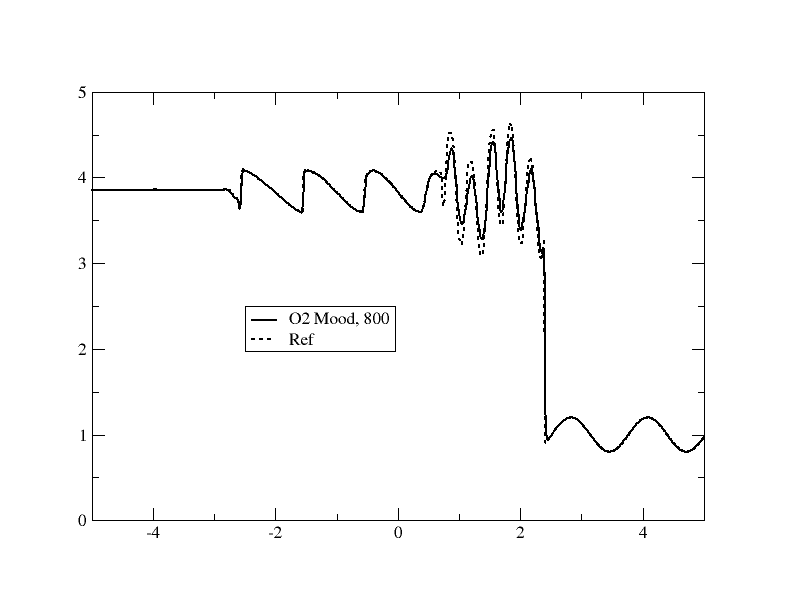}}
		\subfigure[O4Mood]{\includegraphics[width=0.45\textwidth]{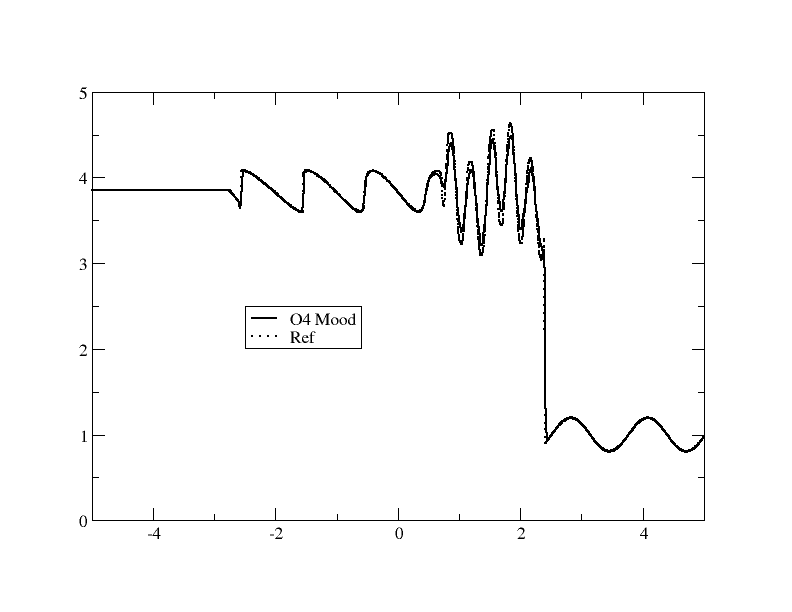}}
		\subfigure[O2MoodNaN]{\includegraphics[width=0.45\textwidth]{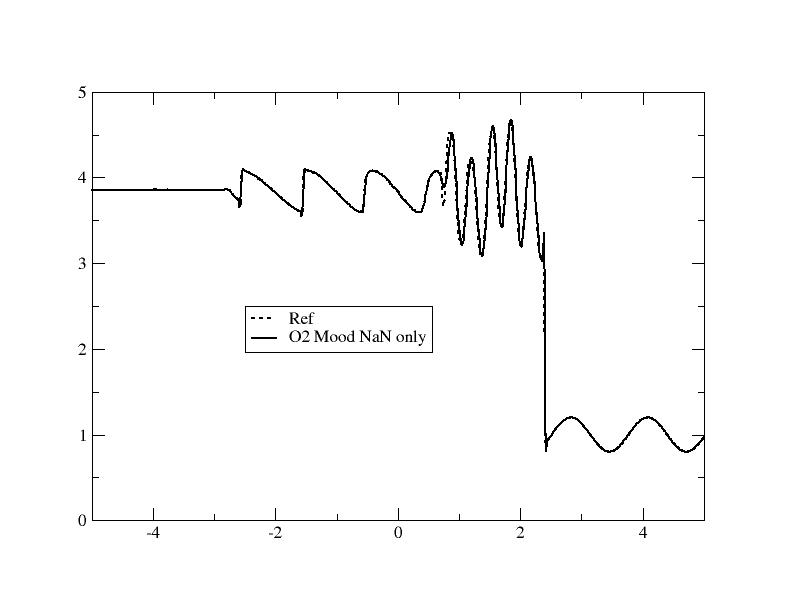}}
		\subfigure[O4MoodNaN]{\includegraphics[width=0.45\textwidth]{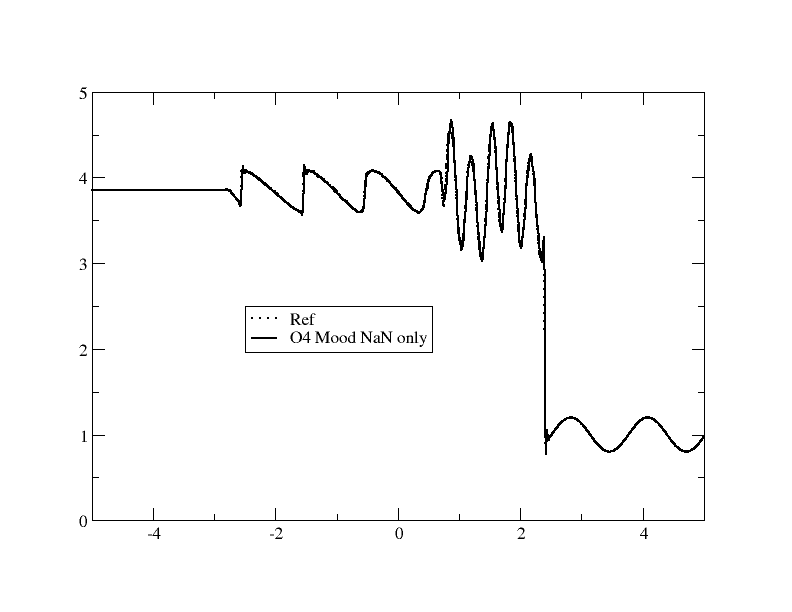}}
		\caption{ The reference solution is plotted with a dotted line. Comparison of various strategies of MOOD, for 800 points.}\label{shu:1}
	\end{center}
\end{figure}

\begin{figure}[h]
	\begin{center}
		\subfigure[O2 800]{\includegraphics[width=0.45\textwidth]{shuRefO2NaN_800.png}}
		\subfigure[O4 800]{\includegraphics[width=0.45\textwidth]{shuRefO4NaN_800.png}}
		\subfigure[O2 400]{\includegraphics[width=0.45\textwidth]{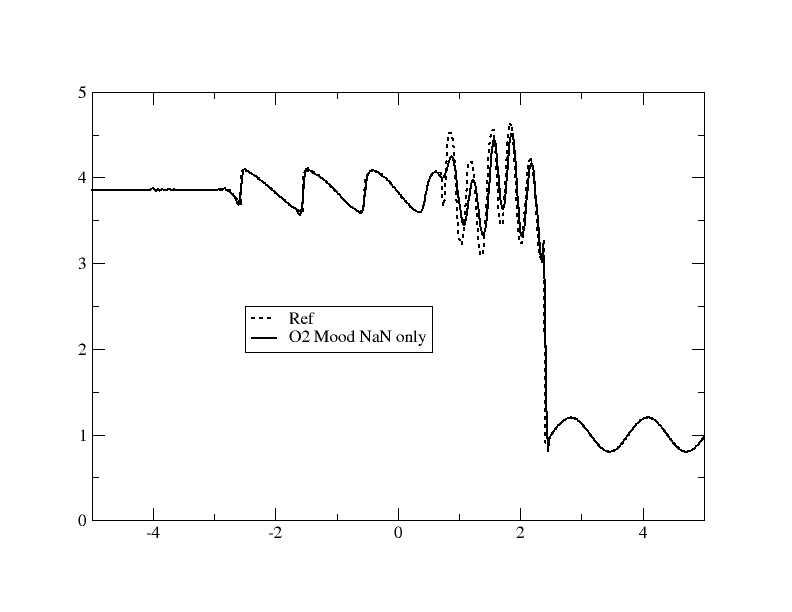}}
		\subfigure[O4 400]{\includegraphics[width=0.45\textwidth]{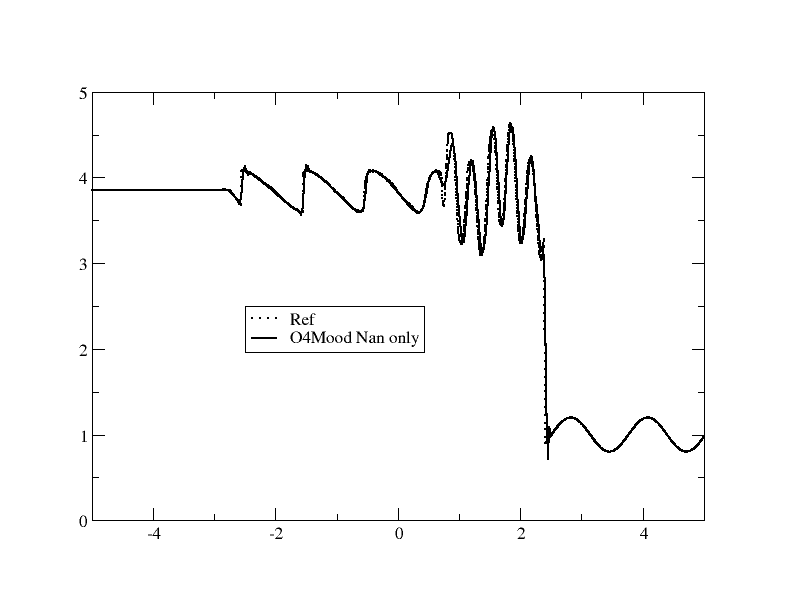}}
		\subfigure[O2 200]{\includegraphics[width=0.45\textwidth]{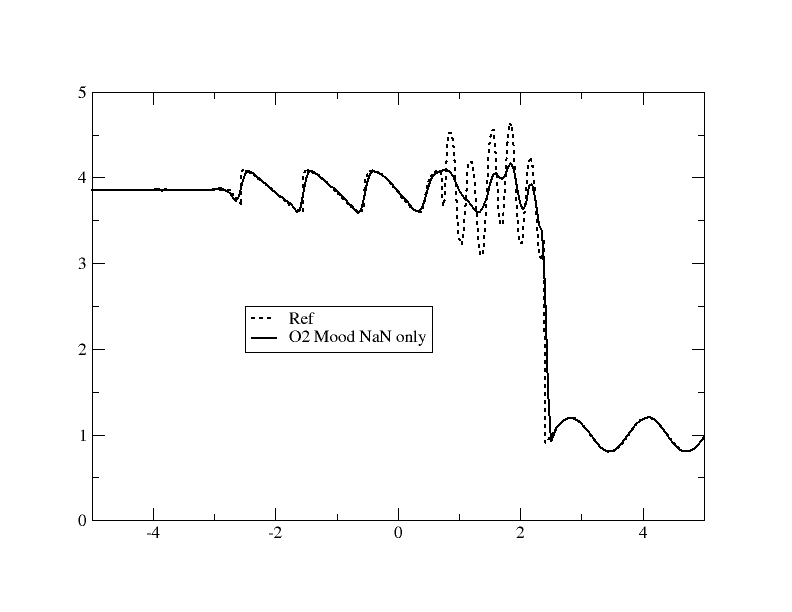}}
		\subfigure[O4 200]{\includegraphics[width=0.45\textwidth]{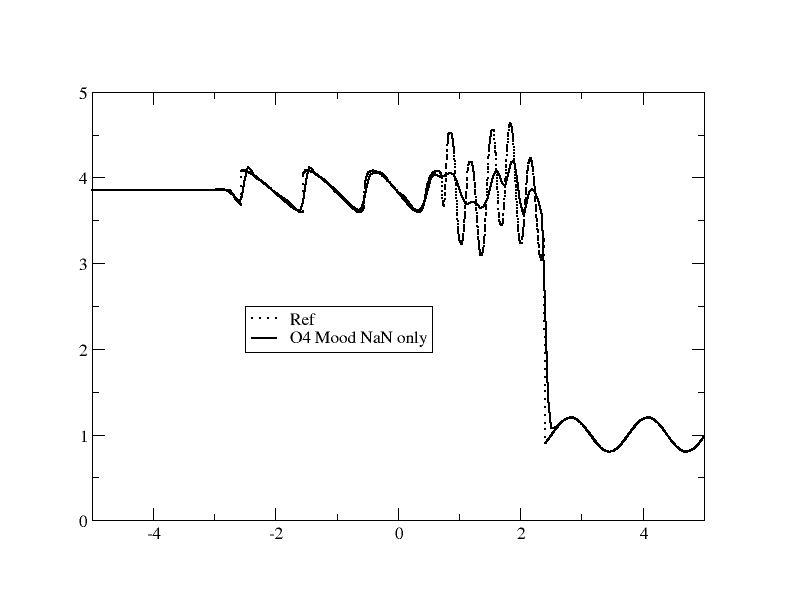}}
		\caption{ O2 and O4 MOOD solutions with control of NaN only, for 200, 400 and 800 mesh points.}\label{shu:2}
	\end{center}
\end{figure}

The conditions of the Shu-Osher test are
$$(\rho, u, p)=\left \{\begin{array}{ll}
	(3.857143, 2.629369, 10.3333333) &\text{if } x<-4,\\
	(1+0.2\sin(5x), 0, 1) &\text{else.}
\end{array}\right .
$$
on the domain $[-5,5]$ and the final time of the problem is $T=1.8$.
The reference solution is obtained with 10.000 points and 4th order limited scheme. It is difficult to see any modification in the solution if we use more grid points, this is why we consider this solution as the reference solution. We display only the solution with the MOOD stabilisation technique, however, we have tried two different strategies. The figures labeled as OXMood, where X=2 or 4, use the full strategy of section \ref{stabilisation}. The physical variables are the density and the pressure, nothing is tested on the velocity. In the figures labeled as OXMoodNaN, with X=2 or 4, we only check if the solution lies in the \revD{invariance domain, i.e., density and pressure stay positive and we do not encounter NaN values. In figure \ref{shu:1} we plot the results for 800 points. In figure \ref{shu:2} we compare results for 200, $400$ and $800$ mesh points only for the MoodNaN strategy.}
From figure \ref{shu:2}, we see that with $800$ points, there is hardly no difference with between the O4MoodNaN solution and the reference one.
%
\section{Conclusion}
In this paper, simplifying a method described in \cite{Torlo}, we show how to construct a class of kinetic numerical methods that can run at least at CFL 1. They can handle in a simple manner hyperbolic problems, and in particular compressible fluid mechanics one. These methods are always locally conservative and thus can handle correctly discontinuities. We have described a rather simple stabilisation mechanism which can be further improved or changed: it is not really the core of the proposed method. Our methodology can be arbitrarily high order and can use CFL number larger or equal to unity on regular Cartesian meshes. Extension to the multidimensional case will be the topic of future works. In particular, our implementation of these methods indicates that they can be potentially very fast. The parallelisation should also be straightforward. This, however, has to be confirmed in several spatial dimensions.

\section*{Acknowledgment}
R.A. would like to thanks Prof. Li-Shi Luo (Old Dominion, USA) and Professor Sagaut (Aix-Marseilles University, France) for their encouragement.
Finally, he would like to thank Professors Lallemand  and D'Humi\`eres for the discussion we had  during his PhD studies. This paper shows that this has never been forgotten. I take this opportunity to also thanks Prof. Chi-Wang Shu for very helpfull discussions about time stepping.
D.T. has been partially funded by  ITN ModCompShock project funded by the European Union's Horizon 2020 research and innovation program under the Marie Sklodowska-Curie grant agreement No 642768.

\appendix
\section{Another iteration technique for the fourth order in time case}
The iteration \eqref{dec} solved by \eqref{solDec} reads (when we have no source term or after the application of the projector operator):
\begin{equation*}
	\begin{split}
		v_1^{(p+1)}=u_i^n-\lambda\big ( \theta_0^1\delta u^n_i+\theta_1^1\delta v_1^{(p)}+\theta_2^1\delta v_2^{(p)}\big )\\
		v_2^{(p+1)}=u_i^n-\lambda\big ( \theta_0^2\delta u^n_i+\theta_1^2\delta v_1^{(p)}+\theta_2^2\delta v_2^{(p)}\big )\\
	\end{split}
\end{equation*}
and after the application of the Fourier transform, we have
\begin{equation*}
	\hv^{(p+1)}=\hu^n e-\lambda g \theta \bv^{(p)}, \quad e=\begin{pmatrix} 1-\theta_0^1 \lambda g\\1-\theta_0^2 \lambda g\end{pmatrix}, \quad
	\theta=\begin{pmatrix} \theta_1^1 &\theta_2^1\\
		\theta_1^2 &\theta_2^2
	\end{pmatrix}
\end{equation*}
so that the amplification factor satisfies
$$G^{(p+1)}=e- \lambda g \theta G^{(p)}.$$
This can be seen as the Jacobi iteration for solving the system 
$$(\text{Id}+\lambda g \theta\big ) G=e.$$
By analogy, we can define a Gauss-Seidel iteration by 
\begin{equation*}
	\begin{split}
		v_1^{(p+1)}=u_i^n-\lambda\big ( \theta_0^1\delta u^n_i+\theta_1^1\delta v_1^{(p)}+\theta_2^1\delta v_2^{(p)}\big )\\
		v_2^{(p+1)}=u_i^n-\lambda\big ( \theta_0^2\delta u^n_i+\theta_1^2\delta v_1^{(p+1)}+\theta_2^2\delta v_2^{(p)}\big )\\
	\end{split}
\end{equation*}
whose Fourier transform is
\begin{equation*}
	\hv^{(p+1)}=\hu^n e-\lambda g \Theta_1 \hv^{(p+1)} -\lambda g \Theta_2\hv^{(p)}, \quad \Theta_1=\begin{pmatrix} 0 & 0\\
		\theta_1^2 &0
	\end{pmatrix}, \quad \Theta_2=\begin{pmatrix} \theta_1^1 &\theta_2^1\\
		0&\theta_2^2
	\end{pmatrix},
\end{equation*}
and hence
$$(\text{Id}+\lambda g\Theta_1)G^{(p+1)}=e-\lambda g\Theta_2 G^{(p)}.$$
In both cases, $G^{(0)}=e$.

We can study the stability of the Gauss--Seidel iteration, and we recall the results of Jacobi's for comparison. Denoting by $g_1$ (resp. $g_2$, $g_{4,1}$, $g_{4,2}$) the Fourier symbol of the operators $a\delta_1$ (resp. $a\delta_2$, $a\delta_4^1$, $a\delta_4^2$), we get the results of table \ref{table:appendix}.
\begin{table}[h]
	\begin{center}
		\begin{tabular}{|c|ccccc|}
			\hline
			Iterations & 1  & 2&3&4&5\\
			\hline
			Symbol &\multicolumn{5}{c|}{Gauss Seidel}\\
			\hline
			$g_1$& 1.5&$1.276906714$&$1.167201858$&$1.197067146$,&$1.152628955$\\
			$g_2$ & 0 & $\geq 1.65$ & $\geq 1.47$ & $\geq 1.435$ & $\geq 1.55$\\
			$g_4^1$\footnote{always 1 for $x=\pi$}&0 & $\geq 0.926$ & $\geq 1.775$ & 0 &0  \\
			$g_4^2$ & 0 & $\geq 0.917$ & $0.8754013933$& $\geq 0.89$ &$\geq 0.86$
			\\
			\hline
			Symbol & \multicolumn{5}{c|}{Jacobi}\\
			\hline
			$g_1$& $1$ & $1$ & $1.256372663$& $1.392646782$ & $1.774161172$\\
			$g_2$&$0$ & $\geq 0.87$ & $\geq 1.625$ & $\geq 1.744$ & $\geq 2.06$\\
			$g_4^1$ \footnote{always 1 for $x=\pi$}&$0$ &$0$ &$\geq 1.25$ & $\geq 2.06$ & $\geq 2.52$\\
			$g_4^2$ &$0$ &$0$ & $\geq 0.905$ & $\geq 1.044$& $\geq 1.321$\\
			\hline
		\end{tabular}
	\end{center}
	\caption{\label{table:appendix} CFL number for stability of the DeC iterations given by Gauss-Seidel and Jacobi methods. $0$ means that the scheme is unconditionally unstable. $x$ means that the scheme is stable up to CFL $x$, $\geq x$ means that the scheme is stable for at least CFL x (and slightly above).}
\end{table}
\begin{remark}[A few remarks about table \ref{table:appendix}.]
	$ $
	\begin{itemize}
		\item For $g_4^1$, the amplification factor is always equal to $1$ when $x=\pi$, and strictly below $1$ under the condition stated above.
		\item For $g_1$ and 3 iterations, the CFL condition can be  computed exactly. It is $\frac{1}{2}\sqrt [3]{4+\sqrt {17}}-\frac{1}{2}\,{\frac {1}{\sqrt [3]{4+\sqrt {17}}}}
		+\frac{1}{2}\approx1.256372663$.
	\end{itemize}
\end{remark}

From this results, we see that there is no fundamental reason to prefer Gauss-Seidel iteration to the Jacobi one; the coding of the Gauss-Seidel is also  slightly more involved. However, this conclusion holds true only for the schemes we have considered here, and might not be true for others.

\bibliographystyle{unsrt}
\bibliography{biblio}
\end{document}